\documentclass[%
times
]{fldauth_noWiley}

\usepackage{graphicx}
\usepackage{amsmath}
\usepackage{amsfonts}
\usepackage{listings}
\usepackage{subfig}
\usepackage{./csml}

\usepackage[pdftex,colorlinks=true,allcolors=blue]{hyperref}

\def\tagform@#1{\maketag@@@{\bfseries(\ignorespaces#1\unskip\@@italiccorr)}}
\renewcommand{\eqref}[1]{\textup{{\normalfont(\ref{#1}}\normalfont)}}
\newcommand{\Mach}{\operatorname{\mathit{M\kern-.15em a}}}
\newcommand{\Reyn}{\operatorname{\mathit{R\kern-.04em e}}}
\newcommand{\Pran}{\operatorname{\mathit{P\kern-.03em r}}}
\newcommand{\Knud}{\operatorname{\mathit{K\kern-.2em n}}}
\newcommand{\CFL}{\operatorname{\mathit{C\kern-.2em F \kern-.2em L}}}

\graphicspath{{./figs/}}

\begin{document}

\runningheads{H Xiao ET AL.}{Immersed discontinuous galerkin method on unstructured meshes}

\title{An Immersed Discontinuous Galerkin Method for Compressible Navier-Stokes Equations on Unstructured Meshes}

\author{Hong Xiao\affil{1,2}, Eky Febrianto\affil{1}, Qiaoling Zhang\affil{1}, Fehmi Cirak\affil{1} \corrauth}

\address{\centering{%
        \affil{1}Department of Engineering, University of Cambridge, Cambridge, CB2 1PZ, UK \\ %
        \affil{2}School of Power and Energy, Northwestern Polytechnical University, Xi'an, 710072, China}
    }

\corraddr{f.cirak@eng.cam.ac.uk}

\begin{abstract}
We introduce an immersed high-order discontinuous Galerkin method for solving the compressible Navier-Stokes equations on non-boundary-fitted meshes. The flow equations are discretised with a mixed discontinuous Galerkin formulation and are advanced in time with an explicit time marching scheme. The discretisation meshes may contain simplicial (triangular or tetrahedral) elements of different sizes and need not be structured. On the discretisation mesh the fluid domain boundary is represented with an implicit signed distance function. The cut-elements partially covered by the solid domain are integrated after tessellation with the marching triangle or tetrahedra algorithms. Two alternative techniques are introduced to overcome the excessive stable time step restrictions imposed by cut-elements. In the first approach the cut-basis functions are replaced with the extrapolated basis functions from the nearest largest element. In the second approach the cut-basis functions are simply scaled  proportionally to the fraction of the cut-element covered by the solid.  To achieve high-order accuracy additional nodes are introduced on the element faces abutting the solid boundary. Subsequently, the faces are curved by projecting the introduced nodes to the boundary. The proposed approach is verified and validated with several two- and three-dimensional subsonic and hypersonic low Reynolds number flow applications, including the flow over a cylinder, a space capsule and an aerospace vehicle. 
\end{abstract}



\maketitle

\section{Introduction}
\begin{figure}[t!]
\centering
    \subfloat[]{\includegraphics[width=0.49\textwidth]{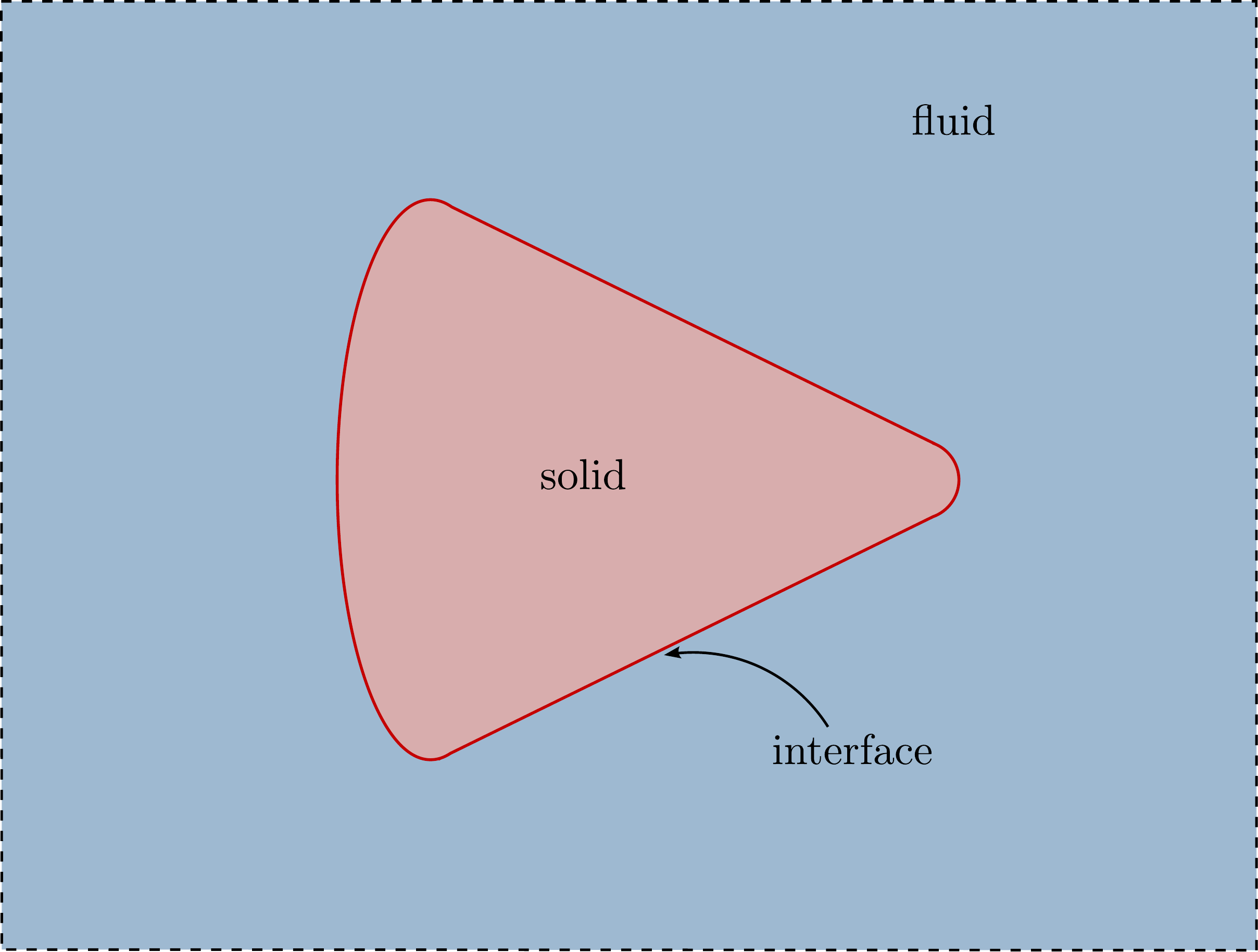}}
    \hfill
    \subfloat[]{\includegraphics[width=0.49\textwidth]{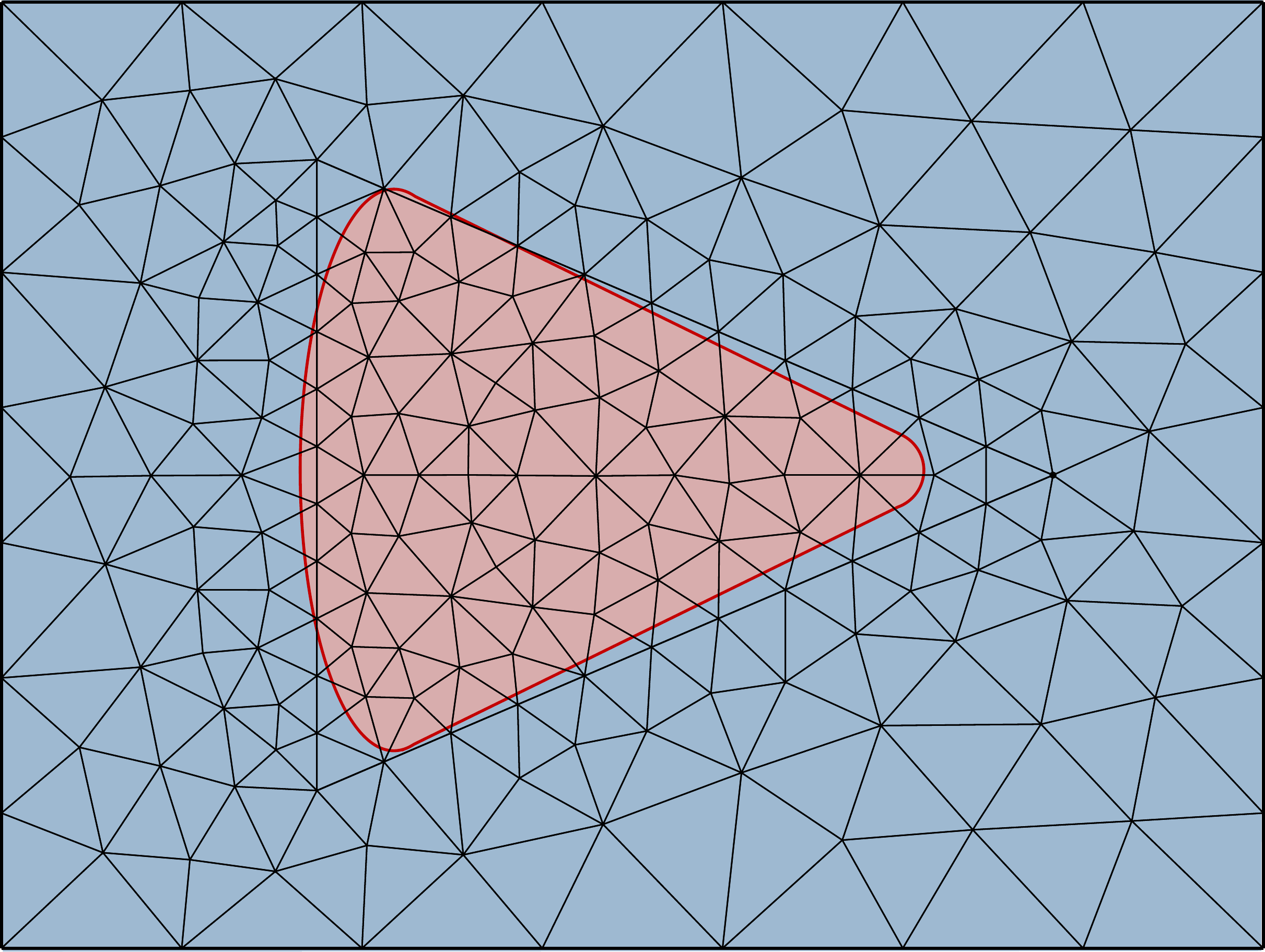}}
  \caption{Immersed discontinuous Galerkin method. (a) A fluid domain (blue) with an immersed solid (red).  (b) Discretisation with a non-boundary-fitted triangular mesh.}
\label{fig:immersedGeneral}
\end{figure}

%
The discontinuous Galerkin method has a number of appealing properties, including accuracy and robustness, when applied to fluid dynamics problems.  In industrial applications with complex geometries a key limitation of computational fluid dynamics using the discontinuous Galerkin method is the need for costly mesh generation. This is common with the finite element and finite volume methods which also require boundary-fitted meshes.  As also identified in the recent NASA CFD Vision 2030 Study, the generation of boundary fitted meshes is  a significant impediment to the implementation of autonomous and reliable computational simulation workflows~\cite{slotnick2014cfd}.  Mesh generation, however, becomes straightforward when the element boundaries need not conform to the domain boundary as in immersed methods, also referred to as embedded or fictitious domain methods, see Figure~\ref{fig:immersedGeneral}. Meshes can still be unstructured and have elements of different sizes in order, for instance, to adapt to critical flow features and boundary layers. In cut-elements intersected by the domain boundary alternative techniques are required to enforce the boundary conditions and evaluate the element integrals. The development of these techniques is greatly aided by the variational structure of the discontinuous Galerkin method. 

Immersed methods, which do not require a boundary-fitted mesh, have a long history in computational fluid dynamics, especially in the context of low order finite difference and finite volume methods. These techniques are often restricted to structured Cartesian grids and can represent the boundaries usually only in a diffuse sense, see e.g.~\cite{mittal:2005, Fedkiw1999, Cirak2007,farhat2008higher}. Critically, their accuracy cannot be systematically improved. Hence immersed methods and high-order accurate schemes appeared until recently to be intrinsically incompatible. Contrary to this perception, lately a number of high-order immersed methods have been introduced mostly for the finite element method, see e.g. ~\cite{duster2008finite, Ruberg2012, ruberg2014fixed, kudela2016smart,lehrenfeld2016high, fries2016higher, main2018shifted}, and a few for the discontinuous Galerkin~\cite{Qin2013, johansson2013high, sun2014adaptive, Muller2017}, finite volume~\cite{ farhat2012fiver, main2017enhanced, huang2018family} and finite difference methods~\cite{Gilmanov:2005aa, choi2007immersed}.  Evidently,  the techniques for the finite element and discontinuous Galerkin methods must share a number of features due to their common origins. Indeed, they all have a method for representing the domain boundaries on the computational mesh, a method for evaluating integrals on elements partially covered by the fluid and a method for enforcing boundary conditions. 

The immersed discontinuous Galerkin method that we have developed builds on earlier work on immersed discontinuous Galerkin and finite element methods. On the computational mesh the domain is described with an implicit signed distance function, similar to that, for instance, in~\cite{Sanches:2011aa,Ruberg2012,Qin2013,Muller2017}.  For integration the cut-elements  are triangulated into simplicial elements using the marching triangle and tetrahedra algorithms~\cite{Gueziec:1995aa, ruberg2014fixed, ruberg2016unstructured}. After the triangulation, element integrals are evaluated with standard Gauss integration. There are alternative techniques available which do not require a cut-element triangulation, but have other limitations~\cite{Muller2017,Schillinger:2011aa,sudhakar2014accurate}.  To achieve high-order accuracy the domain boundaries must be approximated with curved cut-element boundaries of the same polynomial order as the discontinuous Galerkin method basis functions. To this end, we use the curved boundary treatment recently introduced by Zhang~\cite{Zhang2016}. The element-specific domain integrals are evaluated on the simplicial elements with straight boundaries, but the boundary integrals are evaluated on curved boundaries.  In our mixed discontinuous Galerkin implementation, according to Bassi and Rebay~\cite{Bassi1997}, the imposition of the Dirichlet and Neumann boundary conditions on the cut-elements is straightforward. The fluxes over the curved cut-element boundaries are treated in the same way as in the original discontinuous Galerkin method with boundary fitted meshes. 

In most immersed methods numerical accuracy and stability are compromised when only an exceedingly small portion of a cut-element is covered by the fluid. When the discretised equations are solved with an explicit time marching scheme, as in our work, this leads to severe time step size restrictions. In earlier work on immersed discontinuous Galerkin methods element agglomeration strategies have  often been used to deal with this problem. Agglomeration techniques are, however, difficult to implement on unstructured meshes and especially in three dimensions. Borrowing from immersed finite elements, we consider in our work two alternative techniques  to stabilise the small cut-elements. In the first approach the cut-basis functions are replaced with the extrapolated basis functions from the nearest largest element. This is a consequent extension of the extrapolation technique introduced earlier in connection with B-spline basis functions~\cite{Hollig:2001aa, Ruberg2012, ruberg2014fixed}. In the second approach the basis functions of the cut-element are scaled proportionally to the  fraction of the cut-element  covered by the solid. This scaling is motivated by a stabilisation technique proposed in the work of M\"{o}\ss{}ner and Reif~\cite{Mosner2008} for B-spline basis functions. As our numerical computations confirm, both approaches provide an easy to implement and effective means of dealing with cut-element stability issues. Alternatively, as recently proposed, the adaptation of the ghost penalty cut-element stabilisation technique known from immersed, or cut, finite elements to the discontinuous Galerkin setting is also feasible~\cite{gurkan2019stabilized}.

The  paper is organised as follows. In Section 2  the mixed modal discontinuous Galerkin method for compressible Navier-Stokes equation is reviewed. Subsequently, in Section 3, the developed high-order immersed discontinuous Galerkin method is introduced. In Section 4 the application of the immersed discontinuous Galerkin method to problems in subsonic and hypersonic flow regimes is demonstrated. The stability of the two proposed cut-element stabilisation methods is numerically evaluated and the results obtained are compared with those in literature and experiments.  Finally, conclusions and  perspectives are drawn in Section 5.

\section{Mixed modal discontinuous Galerkin method for compressible Navier-Stokes }
%
\subsection{Governing equations}

The non-dimensional compressible Navier-Stokes equations, according to, e.g., Canuto et al.~\cite{canuto2007spectral}, read
\begin{equation}
\frac{\partial{\vec{U}}}{\partial{t}} + \nabla\cdot\vec{F}_{\text{inv}}(\vec{U})
-\nabla\cdot\vec{F}_{\text{vis}}(\vec{U},\nabla\vec{U})= \vec{0} \, 
\label{eq:govEqn}
\end{equation}
with the conservation vector 
\begin{equation}
\renewcommand{\arraystretch}{1.3}
\vec{U}  = \begin{pmatrix}
  \rho \\
  \rho\vec{u} \\
  \rho E
 \end{pmatrix}
\end{equation}
and the inviscid and viscous flux vectors
\begin{subequations}
\begin{align}
\vec{F}_{\text{inv}}(\vec{U}) & =
\renewcommand{\arraystretch}{1.3}
\begin{pmatrix}
  \rho\vec{u} \\
  \rho\vec{u}\otimes\vec{u}+\frac{1}{\gamma \Mach^2}p\vec{I} \\
  \big(\rho E+\frac{1}{\gamma \Mach^2}p\big)\vec{u}
\end{pmatrix} 
\\
\vec{F}_{\text{vis}}(\vec{U},\nabla\vec{U}) &= \frac{1}\Reyn
\renewcommand{\arraystretch}{1.3}
\begin{pmatrix}
  0 \\
  \vec{\Pi} \\
  \vec{\Pi} \vec{u}-\frac{1}{(\gamma-1) \Mach^2Pr}\vec{q}
\end{pmatrix} \, , 
\label{eq:ConservFlux}
\end{align}
\end{subequations}
where $\rho$ is the mass density, $\vec{u}$ the velocity vector, $p$ the pressure, $\gamma$ the specific heat ratio, $\vec I$ the identity tensor, $E$ the total energy density, $\vec{\Pi}$ the deviatoric stress tensor and $\vec{q}$ the heat flux vector. 

The dimensionless Mach number $\Mach$, Reynolds number $\Reyn$ and Prandtl number $\Pran$ are defined as
\begin{eqnarray*}
\Mach = \frac{u_\text{ref}}{(\gamma RT_{\text{ref}})^{1/2}}, \quad \Reyn =\frac{\rho_\text{ref} u_\text{ref} L}{\eta_\text{ref}}, \quad \text{and} \quad \Pran = \frac{C_{p_\text{ref}}\eta_\text{ref}}{k_\text{ref}}.
\label{eq:DimensionlessVar}
\end{eqnarray*}
Here, the reference state variables are the velocity $u_\text{ref}$, temperature   $T_\text{ref}$, density $\rho_\text{ref}$, viscosity $\eta_\text{ref}$, heat capacity  at constant pressure $C_{p_\text{ref}}$ and thermal conductivity $k_\text{ref}$. Furthermore, $R$ is the universal gas constant and  $L$ the characteristic length of the solid. For the viscosity parameter $\eta$ and thermal conductivity $k$, a  power law model in non-dimensional form is assumed
\begin{equation}
\eta=k=T^s \, ,
\label{eq:viscosityAndThermalCond}
\end{equation}
where $s$ is a gas constant associated with the inverse power law of gas molecules. In our examples we take $s=0.78$ for diatomic gas molecules~\cite{Myong1999}.  The constitutive equations for the Navier-Stokes equation are assumed to be
%
\begin{equation}
\vec{\Pi} = \eta[\nabla\vec{u} + (\nabla\vec{u})^\trans] - \frac{2}{3}\eta (\nabla \cdot \vec{u}) \vec{I}, \hspace{10pt} \vec{q} =-k\nabla T \, .
\end{equation}
%
Finally, the ideal gas equation of state $p = \rho T$ is used in the non-dimensional form.

\subsection{Review of the mixed discontinuous Galerkin method}
%
A mixed formulation of the discontinuous Galerkin method was proposed by Bassi and Rebay \cite{Bassi1997} for the treatment of the viscous terms in the Navier-Stokes equation. 
To this end, the governing equation~\eqref{eq:govEqn}, which depends on the second order derivatives of the velocity,  is reformulated as a first-order coupled system.  That is, the auxiliary variable $\vec{S}$ is introduced so that the coupled first order system  can be written as
\begin{subequations}
\begin{align}
& \frac{\partial{\vec{U}}}{\partial{t}}+\nabla\cdot\vec{F}_{\text{inv}}(\vec{U})
-\nabla\cdot\vec{F}_{\text{vis}}(\vec{U},\vec{S}) =0, \\
& \vec{S}(\vec{U}) = T^s \nabla\vec{U}.
\end{align} \label{eq:coupledSys}
\end{subequations}

The discretised field variables $\vec{U}_h$ and $\vec{S}_h$ are in each element $ \omega_j $ represented with
\begin{equation} \label{eq:approximationForm}
\vec{U}_h(\boldsymbol{x},t)=\sum_{i=1}^K \vec{U}^i(t)\varphi_i(\boldsymbol{x}), \quad
\vec{S}_h(\boldsymbol{x},t)=\sum_{i=1}^K \vec{S}^i(t)\varphi_i(\boldsymbol{x}) \, ,
\end{equation}
where~$\varphi_i(\boldsymbol{x})$ are the basis functions, and~$ \vec{U}^i $ and~$ \vec{S}^i $ are the corresponding coefficients. The number of basis functions per element $K$ depends on the polynomial degree of approximation, see Section~\ref{sec:basisFunctions}. In this work the same basis functions~$\varphi_i$ are used in each local fluid element~$\omega_j$, even though it is straightforward to relax this restriction. To obtain the discretised weak form of the governing equations~\eqref{eq:coupledSys}, the discretised coupled system is element-wise multiplied with the test functions~$\varphi_i$ and then integrated by parts over each element $\omega_j$ yielding
\begin{subequations}
\begin{align}
& \frac{d}{d t}\int_{\omega_j} \vec{U}_h \varphi_i \D v
-\int_{\omega_j} (  \vec{F}_{\text{inv}} -  \vec{F}_{\text{vis}}) \cdot \nabla \varphi_i  \D v
+\int_{\partial \omega_j} \varphi_i ( \vec{F}_{\text{inv}}  -  \vec{F}_{\text{vis}})  \cdot \vec n \D s = \vec{0} \label{eq:weakFormulationA} \\
& \int_{\omega_j} \vec{S}_h \varphi_i \D v
+\int_{\omega_j} T^s \vec{U}_h \otimes \nabla \varphi_i  \D v
-\int_{\partial \omega_j} T^s \varphi_i \vec{U}_h  \otimes \vec n \D s=\vec{0} \, , 
\label{eq:weakFormulationB}
\end{align}  \label{eq:weakFormulation}
\end{subequations}
where~$\partial \omega_j$ denotes the boundaries of the element $\omega_j$. 
Due to the discontinuous approximation at element boundaries, the flux functions appearing in the boundary integrals of each element are replaced by numerical flux functions. We approximate the inviscid term~$\vec{F}_{\text{inv}} \cdot \vec n $ with the Lax-Friedrichs flux 
\begin{equation} \label{eq:LaxFriedrichs}
\vec{h}_{\text{inv}}(\vec{U}_h^-,\vec{U}_h^+ ; \vec{n})=
\frac{1}{2}\Big[\vec{F}_{\text{inv}}(\vec{U}_h^-)+\vec{F}_{\text{inv}}(\vec{U}_h^+) \Big] \cdot \vec{n}  - \frac{C}{2}(\vec{U}_h^+-\vec{U}_h^-)
\end{equation}
with  
\begin{equation}
C = \max\Big(\big| \vec u^- \cdot \vec n \big|+\frac{c^-}{\Mach},\big| \vec u^+ \cdot \vec n \big|+\frac{c^+}{\Mach}\Big) \, ,
\end{equation}
where $c=T^{1/2}$ is the non-dimensional speed of sound. The superscripts $-$ and $+$ denote the variables belonging to the elements on the left and right of an edge. Note that the appearance of $\Mach$ in the definition of $ C $ is due to the used dimensionless formulation. For other more advanced flux approximation schemes see, e.g., the work of Toro~\cite{toro2013riemann}. As a flux limiter, in case of linear basis functions the technique presented in the work of Cockburn et al.~\cite{Cockburn1988} is applied direction-by-direction and in case of quadratic basis functions a coupled slope limiter similar to the one in Le et al.~\cite{le2014triangular} is used.

We approximate the viscous flux~$\vec{F}_{\text{vis}} \cdot \vec n $ with the central flux 
\begin{equation}
\vec{h}_{\text{vis}}(\vec{U}_h^-,\vec{S}_h^-,\vec{U}_h^+,\vec{S}_h^+; \vec{n}) = 
\frac{1}{2}\Big[\vec{F}_{\text{vis}}(\vec{U}_h^-,\vec{S}_h^-)+
\vec{F}_{\text{vis}}(\vec{U}_h^+,\vec{S}_h^+)\Big]\cdot \vec{n},
\label{eq:viscousFlux}
\end{equation}
and the term $T^s\vec{U}_h  \otimes \vec n$ in the auxiliary equation~\eqref{eq:weakFormulationB} with
\begin{equation}
\vec{h}_{\text{aux}}(\vec{U}_h^-,\vec{U}_h^+; \vec{n}) = 
\frac{1}{2}\Big[T^{s-}\vec{U}_h^-+
T^{s+}\vec{U}_h^+\Big] \otimes \vec{n} \, .
\label{eq:auxiliaryFlux}
\end{equation}

All the integrals in the preceding discretised equations are evaluated with Gaussian quadrature. This leads to a semi-discretised system of equations, which can be compactly written as
\begin{subequations} \label{eq:matrixSys}
\begin{align}
\mat{L}\frac{ \D \mat{U}}{\D t} &=\mat{R}_{\mat U}(\mat{U}), \label{eq:matrixSys1} \\
\mat{L} \mat{S} &= \mat{R}_{\mat S}(\mat{U}) \, ,
\label{eq:matrixSys2}
\end{align} 
\end{subequations}
where $\ary L$ is the mass matrix and $ \ary U$ and $\ary S$ are arrays containing the coefficients of the basis functions.  The corresponding residual vectors are denoted with $\mat{R}_{\mat U}$ and $\mat{R}_{\mat S}$. We solve this system using a third-order five-step Runge-Kutta time integration \cite{Cockburn2001}.  The explicit time integration scheme allows us to decouple the solution of~\eqref{eq:matrixSys} in two steps. Firstly, equation~\eqref{eq:matrixSys2} for the auxiliary unknown $\vec S_h$  is solved to compute the derivatives $\nabla \vec U_h$ of the conservative variables. Subsequently, equation~\eqref{eq:matrixSys1} is solved with the already computed  $\nabla \vec U_h$.  As will be discussed in Section \ref{sec:basisFunctions}, we choose orthogonal basis functions so that the mass matrix $\mat{L}$ is readily invertible.  

The time step size $\Delta t$ is chosen according to Cockburn and Shu~\cite{Cockburn2001} with 
\begin{equation}
\Delta t=\min \{ \Delta t_j \}, \quad
\Delta t_j=\frac{1}{2m+1}\frac{ h_j \, \CFL}{|u|+c /\Mach+\eta/(\rho h_j \Reyn)} \, ,
\label{eq:timeSteps}
\end{equation}
where the index $j$ refers to the element $\omega_j$ with the characteristic element size $ h_j $, $m$ is the polynomial degree of the basis functions and $\CFL < 1$ is the Courant-Friedrichs-Lewy (CFL) number. 
%
\subsection{Boundary and initial conditions \label{sec:bouInitial}}
%
The boundary conditions for the compressible Navier-Stokes equations~\eqref{eq:govEqn} depend critically on the Knudsen number $\Knud$, which is the ratio of molecular mean free path $\lambda$  to the characteristic length of the solid $L$ and can be expressed with 
\begin{equation}
	\Knud = \frac{\lambda}{L} = \sqrt{\frac{\gamma \pi}{2} }\frac{\Mach}{\Reyn} \, .
\end{equation}
Especially, for rarefied (large $\lambda$) or microscale gas flows (small $L$) the typical non-slip boundary conditions are not suitable.  To cover all Knudsen number regimes, we consider the Langmuir boundary conditions that are derived from the Langmuir adsorption isotherm at the gas and boundary surface interface  \cite{Choi2005, myong2004gaseous, Mahdavi2014}. 

The fraction of the boundary surface covered by gas molecules is given by the ratio~$\alpha$, with \mbox{$0 \le \alpha \le 1$,} and can be expressed in dimensional form as 
\begin{equation}
\alpha=
\renewcommand{\arraystretch}{1.3}
\begin{cases}
\dfrac{\beta p}{1+\beta p} & \text{for  monatomic} \\[10pt]
\dfrac{\sqrt{\beta p}}{1+\sqrt{\beta p}} & \text{for  diatomic},
\end{cases}
\label{eq:coverageFrac}
\end{equation}
where $p$ is the pressure at the boundary and $\beta$ is an equilibrium constant that describes the interaction between the boundary and gas molecules. The equilibrium constant may be derived as  
\begin{equation}
\beta=\sqrt{\frac{\pi}{32}}\frac{\pi}{c_i^2}\frac{T_r}{T_w} \exp \bigg(\frac{D_e}{R T_w}\bigg)\frac{1}{p_r \Knud},
\label{eq:beta}
\end{equation}
where $c_i$ is the exponent of the inverse power law for the gas particle interaction potential,  $p_r$ the reference pressure, $T_r$ the reference temperature, $T_w$ the wall temperature, $D_e$ the heat of adsorption and $R$ the universal gas constant. The two constitutive constants in this model for air and nitrogen are $c_i = 1.1908$ and  $D_e=5255  \; J/mol$. 

The boundary conditions for the velocity and temperature are given in dependence of the ratio~$\alpha$ as follows  
\begin{subequations} \label{eq:slipJumpBC}
\begin{align}
\vec{u} & =\alpha\vec{u}_w+(1-\alpha)\vec{u}_g \label{eq:slipJumpBC1}\\
T &=\alpha T_w+(1-\alpha) T_g, \label{eq:slipJumpBC2}
\end{align} 
\end{subequations}
where $\vec{u}_w$ and $\vec{u}_g$ denote  wall velocity and far-field gas velocity respectively; and, $T_g$ is far-field gas temperature. As can be inferred from \eqref{eq:coverageFrac} and \eqref{eq:beta}, the ratio~$\alpha$ approaches  $ \alpha = 1 $ with the decrease in Knudsen number $\Knud$ so that the boundary condition according \eqref{eq:slipJumpBC} becomes no-slip.

In the mixed discontinuous Galerkin method with explicit time marching  it is straightforward to apply the wall boundary conditions~\cite{Bassi1997}. The inviscid and auxiliary fluxes in \eqref{eq:LaxFriedrichs} and \eqref{eq:auxiliaryFlux} are computed with the  slip velocity~\eqref{eq:slipJumpBC1}, the temperature~\eqref{eq:slipJumpBC2} and other variables taken from inside the fluid, that is, 
\begin{subequations}
\begin{align}
\vec{h}_{\text{inv}}(\vec{U}^-,\vec{U}^+; \vec{n})\Big|_{\text{w}} &=
\vec{F}_{\text{inv}}(\vec{U } |_{\text{w}})\cdot \vec{n} \, ,
\\
\vec{h}_{\text{aux}}(\vec{U}^-,\vec{U}^+;\vec{n})\Big|_{\text{w}} &=
T^s |_{\text{w}}\vec{U} |_{\text{w}}\otimes \vec{n} \, .
\label{eq:wallBCeval}
\end{align}
\end{subequations}
The viscous flux in  \eqref{eq:viscousFlux} is evaluated using the value of the unknown auxiliary variable $\vec{S}$,
\begin{equation}
\vec{h}_{\text{vis}}(\vec{U}^-,\vec{S}^-,\vec{U}^+,\vec{S}^+; \vec{n})\Big|_{\text{w}}=
\vec{F}_{\text{vis}}(\vec{U}_{\text{w}},\vec{S}_{\text{w}})\cdot \vec{n}.
\label{eq:vicousBCeval}
\end{equation}
The initial condition for the conservation variables throughout the domain is specified using values at the far-field. The coefficients of the basis functions in the interpolation equation~\eqref{eq:approximationForm} are initialised as follows
\begin{equation}
\vec U^1 \Big|_{t=0} = \vec{U}_{\infty}, \qquad  \vec U^i\Big|_{t=0} = 0 \quad \text{for } i = 2,...,K \, .
\label{eq:initialConditionU}
\end{equation}
Note that for the considered modal basis functions we have $\varphi_1(\vec x)=1$. The components of auxiliary unknown are set to zero, 
\begin{eqnarray*}
\vec S^i\Big|_{t=0} = 0  \quad \text{for } i=1,2,...,K.
\label{eq:initialConditionS}
\end{eqnarray*}
%
\subsection{Modal basis functions \label{sec:basisFunctions}}
%
We discretise the two field variables $\vec U$ and $\vec S$, c.f.~\eqref{eq:coupledSys}, using the modal Dubiner basis functions~\cite{Dubiner1991}, which are derived from the Jacobi polynomial
\begin{equation}
P_n^{\alpha,\beta}(x)=\frac{(-1)^n}{2^nn!}(1-x)^{-\alpha}(1+x)^{-\beta}\frac{\D^n}{\D x^n}
\big[(1-x)^{n+\alpha}(1+x)^{n+\beta}\big] 
\label{eq:JacobiPoly}
\end{equation}
with $\{ \alpha, \beta \} \in \mathbb N_0$, not to be confused with the same symbols used in Section~\ref{sec:bouInitial}. 
In addition, the following three functions are defined 
\begin{subequations}
\begin{align}
\Theta_i^a(x) &=P_i^{0,0}(x),\\
\Theta_{ij}^b(x) &=\Big(\frac{1-x}{2}\Big)^iP_j^{(2i+1),0}(x),\\
\Theta_{ijk}^c(x) &=\Big(\frac{1-x}{2}\Big)^{i+j}P_k^{(2i+2j+1),0}(x).
\label{eq:JacobiPolyProds}
\end{align}
\end{subequations}

\begin{figure}
  \centering 
  \includegraphics[width=0.95\textwidth]{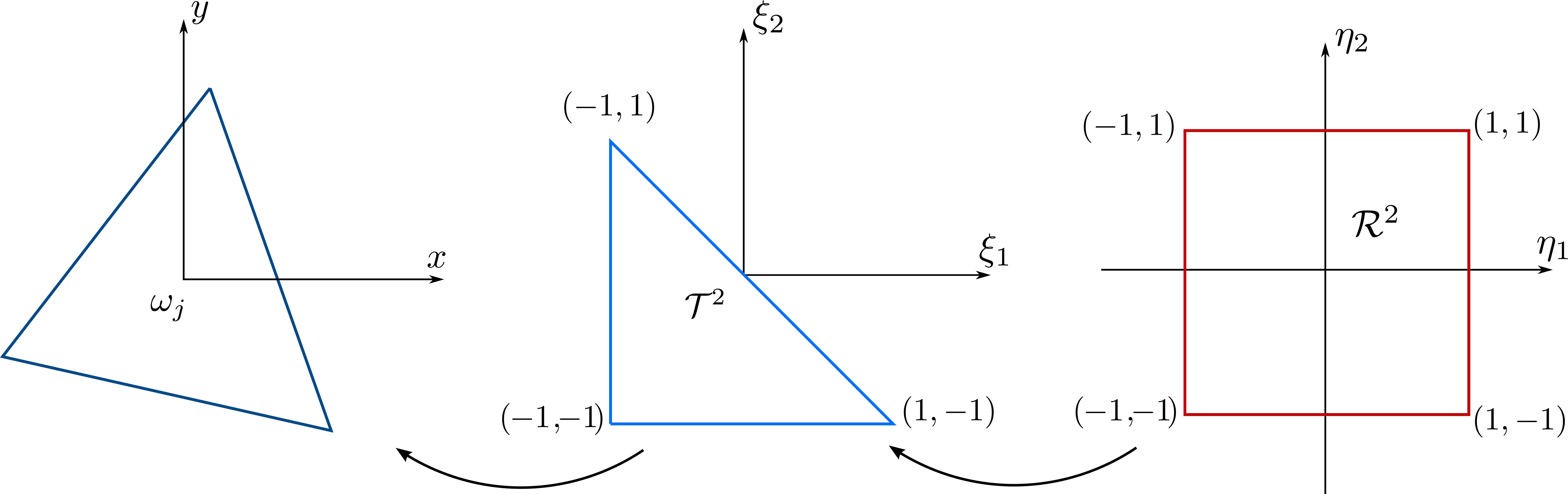}
  \caption{Mapping of the reference square $\mathcal{R}^2$ square to the reference triangle  $\mathcal{T}^2$ and the physical triangle  $\omega_j$.}
  \label{fig:transform2D}
\end{figure}

In two-dimensions, as illustrated in Figure~\ref{fig:transform2D}, we firstly construct the basis in a reference square 
\begin{equation}
\mathcal{R}^2 = \{ (\eta_1,\eta_2) \in \mathbb{R}^{2}; -1 \le \eta_1 \le 1; -1 \le \eta_2 \le 1 \} \, 
\label{eq:quadrilateralT}
\end{equation}
and then map it to the reference triangle 
\begin{equation}
\mathcal{T}^2 = \{ (\xi_1,\xi_2) \in \mathbb{R}^{2}; \xi_1 \ge -1; \xi_2 \ge -1; \xi_1 + \xi_2 \le 0 \} \, .
\label{eq:triangleT}
\end{equation}
To this end, the linear mapping 
\begin{equation}
	 \eta_1 = \frac{ 1 + 2 \xi_1 + \xi_2 }{ 1 - \xi_2 } \text{ \; and \;}  \eta_2 = \xi_2 
\end{equation}
is used, which collapses one of the edges of the square. The orthogonal basis on $ \mathcal{T}^2 $ is given by  
\begin{equation}
\varphi_{(i,j)}( \xi_1, \xi_2)=\Theta_i^a(\eta_1)\cdot\Theta_{ij}^b(\eta_2) \, .
\label{eq:DubinerBasisR2}
\end{equation}

\begin{figure}
  \centering
  \includegraphics[width=0.95\textwidth]{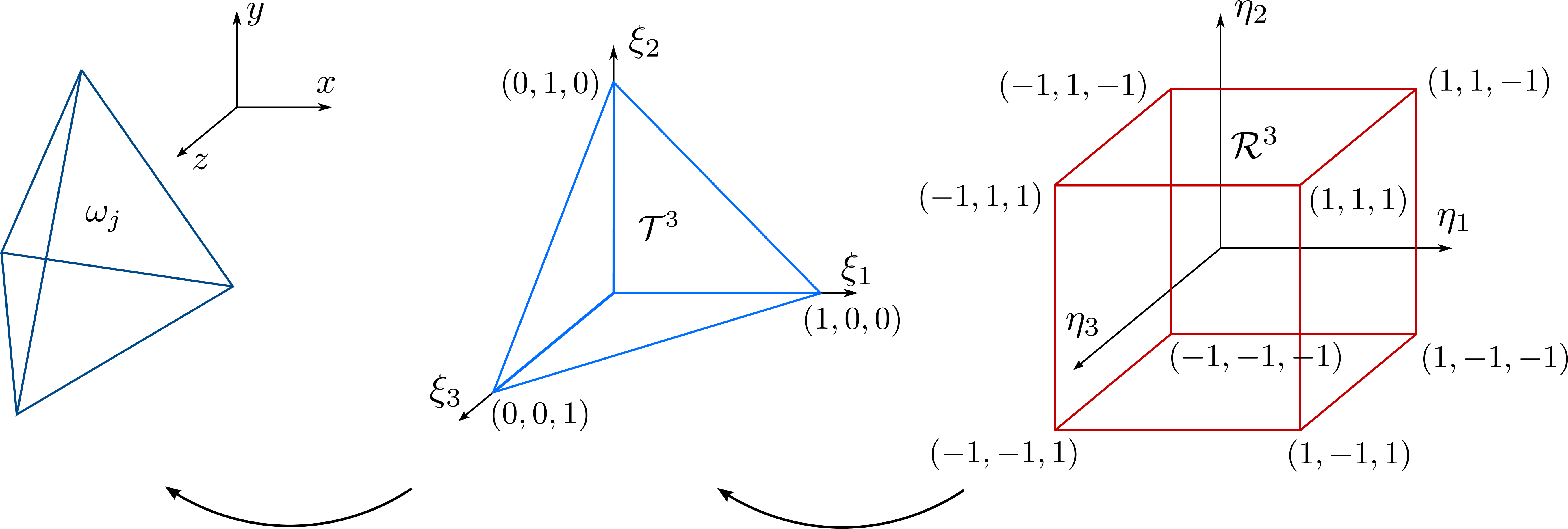}
  \caption{Mapping of the reference cube  $\mathcal{R}^3$ to the reference tetrahedron $\mathcal{T}^3$ and the physical element $\omega_j$. }
  \label{fig:transform3D}
\end{figure}

Similarly, in three-dimensions as shown in Figure~\ref{fig:transform3D} the basis functions are first constructed in a reference hexahedron 
\begin{equation}
\mathcal{R}^3 = \{ (\eta_1, \eta_2, \eta_3) \in \mathbb{R}^{3}; 0 \le \eta_1 \le 1; 0 \le \eta_2 \le 1 ; 0 \le \eta_3 \le 1 \} \,
\label{eq:hexahedralT}
\end{equation}
and then mapped into a reference tetrahedron 
\begin{equation}
\mathcal{T}^{3} = \{ (\xi_1,\xi_2,\xi_3) \in \mathbb{R}^{3}; \xi_1 \ge 0; \xi_2 \ge 0; \xi_3 \ge 0; \xi_1 + \xi_2 + \xi_3 \le 1 \} \, .
\label{eq:tetrahedraT}
\end{equation}
The orthogonal basis functions in the reference tetrahedron are given by 
\begin{equation}
\varphi_{(i,j,k)}(\xi_1, \xi_2, \xi_3) = \Theta_i^a(\eta_1) \cdot \Theta_{ij}^b(\eta_2) \cdot \Theta_{ijk}^c (\eta_3)
\label{eq:DubinerBasisR3}
\end{equation}
with the mappings
\begin{equation}
 \eta_1 = \frac{-1 + 2\xi_1 + \xi_2 + \xi_3}{1 - \xi_2 - \xi_3},  {\; } \eta_2=\frac{-1 + 2\xi_2 + \xi_3}{1-\xi_3} \text{ \; and \;}  \eta_3=-1 + 2\xi_3  \, .
 \end{equation}

%

%

%

%
\section{Immersed discontinuous Galerkin method on non-boundary-fitted meshes}
%

\subsection{Implicit geometry description}
%
To begin with, we assume that the embedded solid domain $\Omega$ is described with a scalar-valued implicit (or, level set) function 
\begin{equation}\label{eq:sgndist1}
   \phi (\vec{x}) = 
   \begin{cases} 
     -   \dist( \vec{x}, \Gamma) \quad & \text{if } \vec x  \in \Omega \\ 
     \phantom{-}0 & \text {if } \vec{x} \in \Gamma \\
     \phantom{-} \dist(\vec{x},\Gamma) & \text{otherwise}\,,
   \end{cases}
\end{equation}
where $\dist(\vec x,\Gamma) = \min_{\vec y \in \Gamma} |\vec x - \vec y|$  is  the shortest distance between the point $\vec x$ and the fluid-solid interface $\Gamma$. By definition, the level set function $\phi(\vec x)$ is negative inside the solid domain $\Omega$ and positive inside the fluid domain. The unit normal to the boundary is given by the normalised gradient of the level set function
\begin{equation} \label{eq:sgndistNorm}
	\vec{n}(\vec x) = \frac{\nabla \phi (\vec x)}{ \| \nabla \phi (\vec x) \|}	 \, ,
\end{equation}
where $\vec x \in \Gamma$ is a point on the boundary.

\begin{figure}
  \centerline{\includegraphics[width=0.55\textwidth]{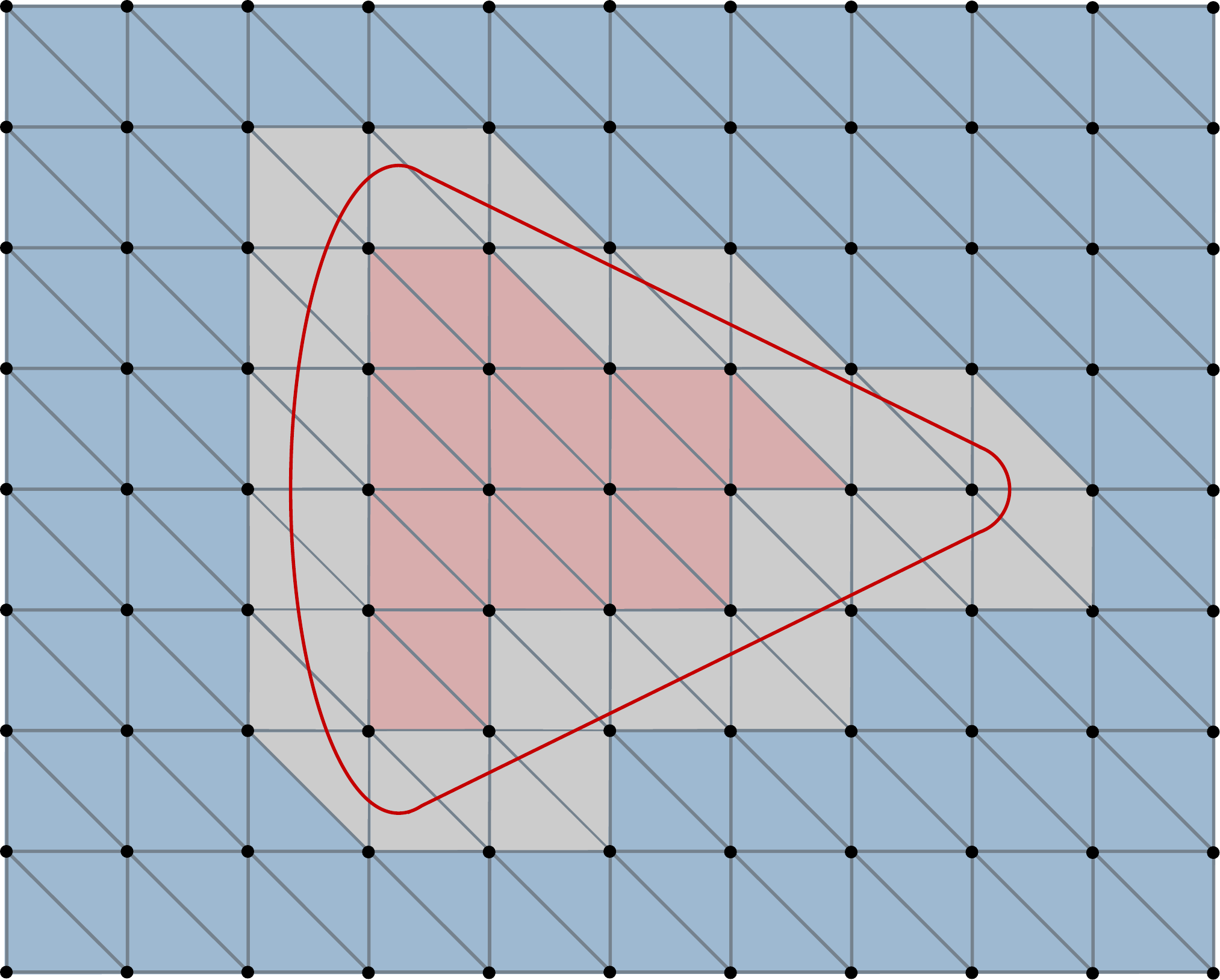}}
  \caption{Discontinuous Galerkin mesh with an embedded solid. The fluid solid interface is outlined in red. The elements are classified as fluid (blue), solid (red) or cut (grey).}
  \label{fig:element categorisation}
\end{figure}

In problems where an implicit function $\phi(\vec x)$ is not available a sufficiently accurate signed distance function needs to be computed. This can be accomplished either with numerical signed distance computations~\cite{jones20063d,Mauch:2003aa,ruberg2014fixed} or analytic implicitisation techniques~\cite{upreti2014algebraic}.
The use of an implicit instead of a parametric representation greatly simplifies all geometric operations. For instance, the classification of elements as solid, fluid or cut involves only the evaluation of the signed distance function $\phi(\vec x)$, see Figure \ref{fig:element categorisation}. In our present implementation elements are classified according to the signed distance function value at their nodes. In  the discontinuous Galerkin computations only the fluid and cut-elements are considered. 

%

%
\subsection{Cut-element  integration \label{sec:cutElemIntegration}}
%
For elements entirely inside the fluid domain the element and boundary integrals appearing in the weak form of the governing equations~\eqref{eq:weakFormulation} are evaluated with Gaussian quadrature, see  Section~\ref{sec:resultDiscussion}. The cut-elements, which are only partly covered by the fluid domain, are first triangulated before evaluating the domain and boundary integrals. Note that this triangulation is only for integration. The field variables are still approximated by the modal basis functions defined over the entire cut-element.

\begin{figure}
  \centerline{\includegraphics[width=0.75\textwidth]{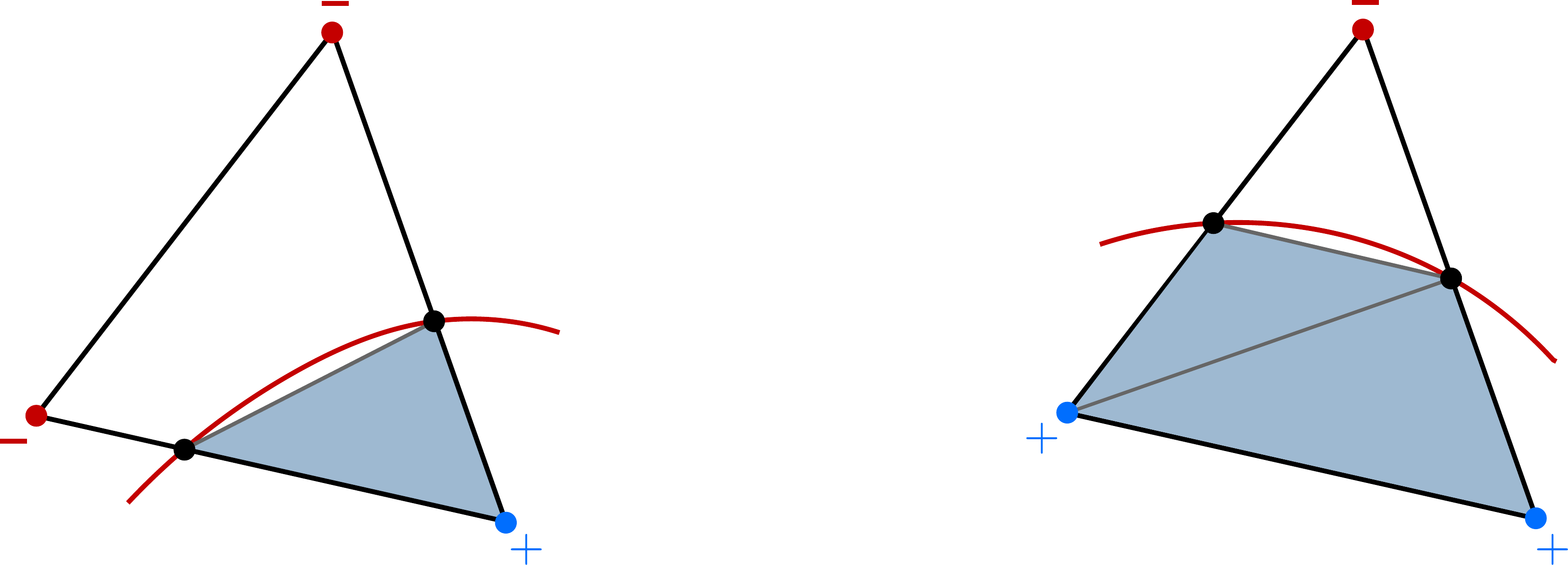}}
  \caption{Two canonical cases of how a triangle may be traversed by the boundary (red line). Red nodes have negative level set values and blue nodes have positive level set values. The appropriate triangulation of the fluid part for both cases is also given in the two figures.}
  \label{fig:marchingTri}
\end{figure}

\begin{figure}[h]
	\centering
	\subfloat[]{\includegraphics[width=0.9\textwidth]{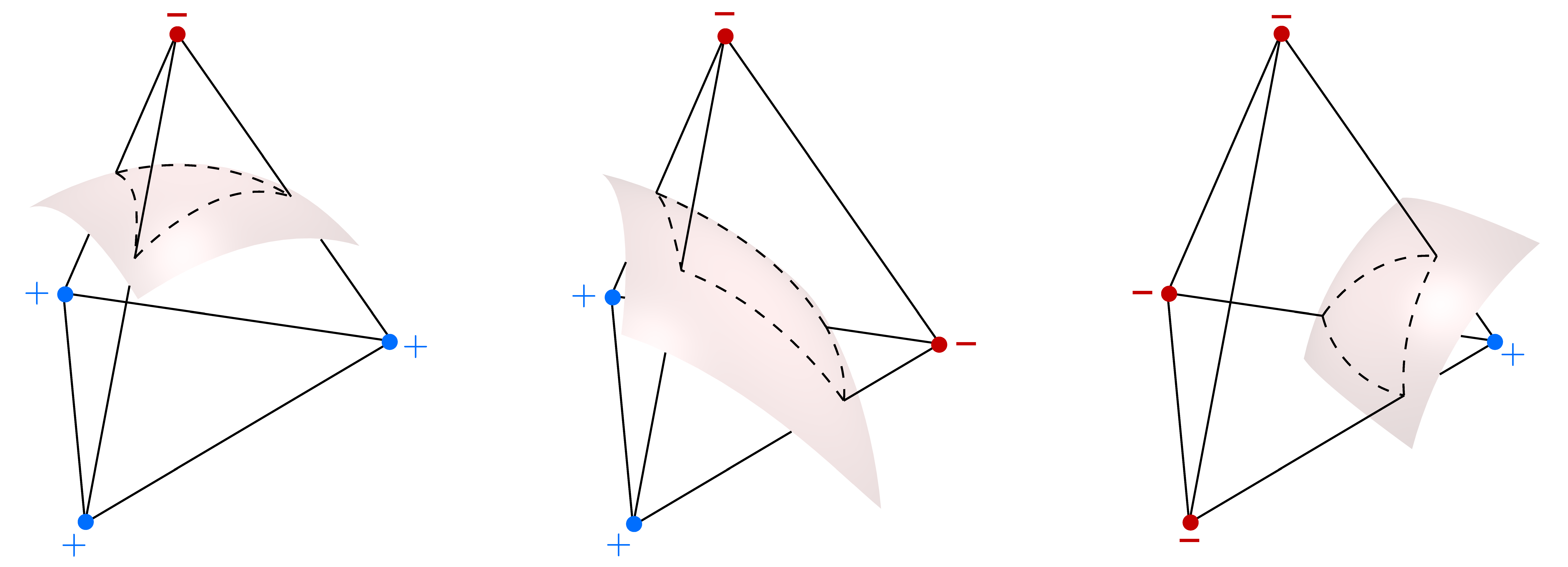}}	
	\vspace{0.01\textwidth}
	\subfloat[]{\includegraphics[width=0.9\textwidth]{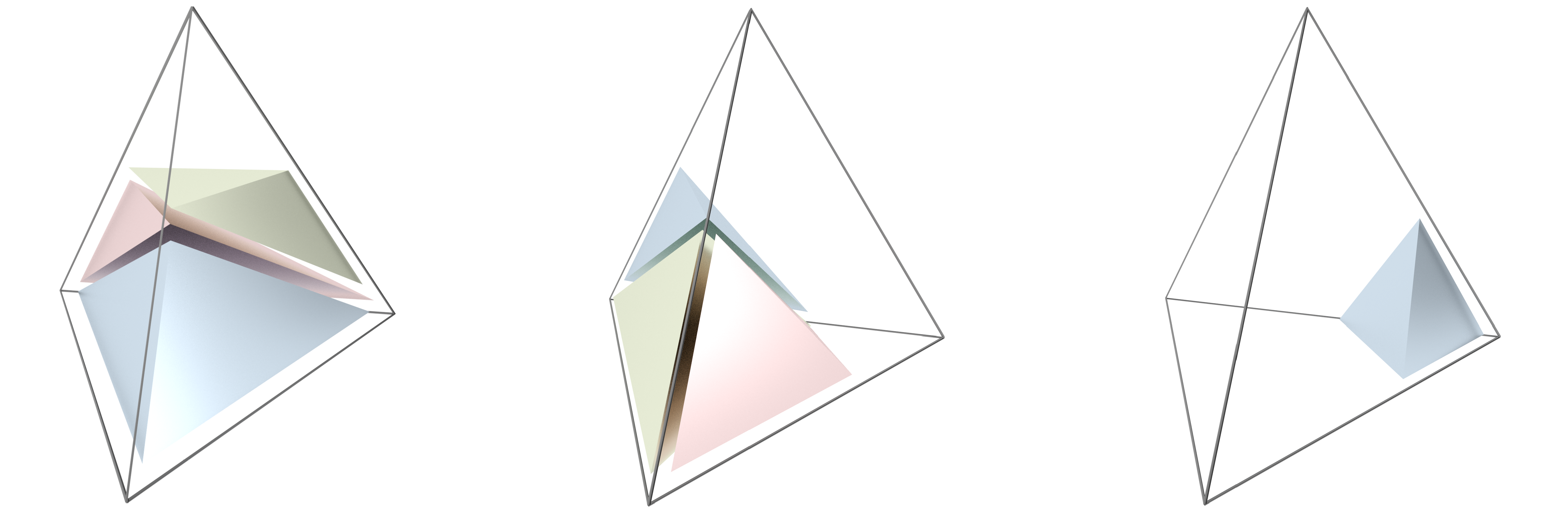}}	
  \caption{(a) Three canonical cases of how a tetrahedron may be traversed by the boundary (red surface). Red nodes have negative level set values and blue nodes have positive level set values. (b) The appropriate tetrahedralisation of the fluid part for the three cases.  \label{fig:marchingTet}}
\end{figure}

In our implementation the cut-elements are triangulated according to R\"uberg and Cirak \cite{ruberg2014fixed} using marching triangle in 2D and marching tetrahedra algorithms in 3D. In 2D there are  $2^3$ combinations for the sign of the level set function at the three nodes of the triangle. Considering symmetry these reduce to the  two canonical cases shown in Figure~\ref{fig:marchingTri}. In 3D, there are $2^4$ combinations for the sign of the level set function, which reduce to the three canonical cases shown in Figure~\ref{fig:marchingTet}. In evaluating the element and flux integrals, we use in each sub-triangle or sub-tetrahedron the same number of Gauss integration points as indicated in Section~\ref{sec:resultDiscussion} for the elements fully covered by the fluid. 
Hence, the total number of integration points for cut-elements is usually higher than for other elements. This does however not lead to a significant time overhead because only very few elements in a mesh are cut-elements. Moreover, while the elements are classified as fluid, solid or cut according to the signed distance function value at their nodes, the intersection of the boundary with the level set function is found by bisectioning.

%

%
\subsection{Cut-element stabilisation}
%
When an exceedingly small portion of a cut-element is covered by the fluid domain, numerical accuracy and stability are usually compromised. The diagonal entries of the  mass matrix $\ary L$ in the system of  equations~\eqref{eq:matrixSys} are proportional to the physically active support size of the corresponding basis functions. Consequently, the condition number of  the matrix $\ary L$ can become arbitrarily large if there are cut-elements with only a negligible part covered by the fluid domain.  This results in significantly smaller stable time step sizes $\Delta t$ when an explicit time marching scheme is used. Although no estimates exist for such cut-elements, it can be argued that the characteristic element size $h_j$ in the estimate~\eqref{eq:timeSteps} has to be chosen in accordance with the physically active support size of each basis function, i.e. the area in 2D or the volume in 3D of the triangulated part of the cut-element.

We propose two alternative stabilisation techniques to overcome the excessive stable time restriction imposed by cut-elements. Both stabilisation techniques can also be interpreted as preconditioning of the  system of discretised discontinuous Galerkin equations~\eqref{eq:matrixSys}. To decide which cut-elements are susceptible to stability issues, the following coverage ratio is computed
\begin{equation}
f_j=\frac{\text{Vol}(\omega_j^F)}{\text{Vol}(\omega_j)} \, ,
\label{eq:fraction}
\end{equation}
where $\text{Vol}(\omega_j^F)$ is the size, area in 2D or volume in 3D, of the triangulated part of the cut-element and the $\text{Vol}(\omega_j)$ is the size of the entire cut-element.  Although, in principle, all basis functions belonging to cut-elements with $f_j < 1$ may lead to stability issues, we process in our implementation for numerical efficiency reasons only cut-elements with $f_j <  0.3$.

\begin{figure}[]
  \centerline{\includegraphics[width=0.4\textwidth]{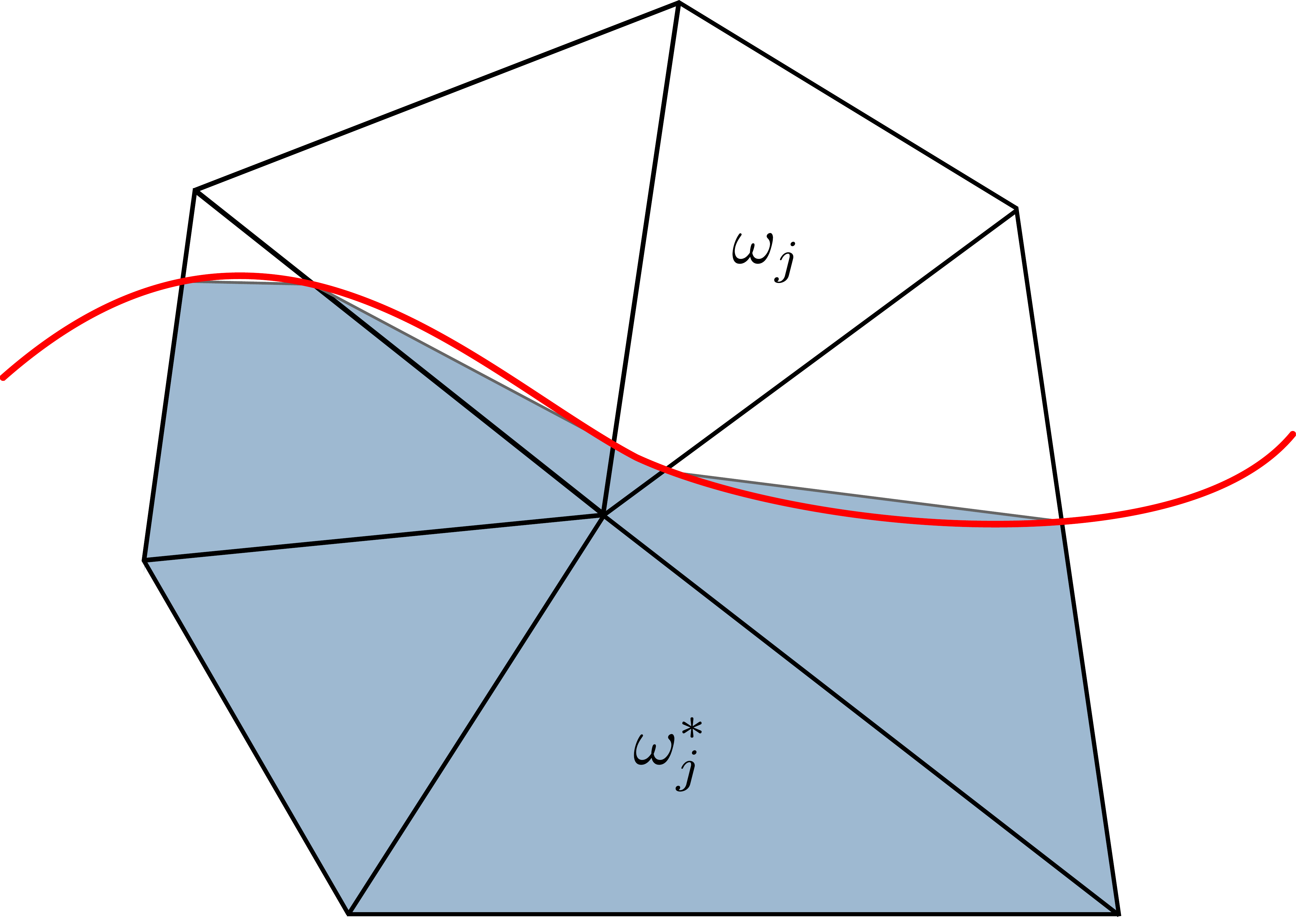}}
  \caption{Cut-element stabilisation by extension. The cut-element $ \omega_j $ with a very small portion covered by the fluid domain is a candidate for stabilisation. The discretised conservation vector~$\vec U_h$ and auxiliary vector~$\vec S_h$ in the cut-element are extrapolated from the basis functions and coefficients} of the largest neighbouring fluid element $\omega^*_j$.  
  \label{fig:basisExtension}
\end{figure}

In the first proposed stabilisation approach, we replace the basis functions in the critical cut-element with the extrapolated basis functions of the neighbouring fluid element of largest size, see Figure~\ref{fig:basisExtension}. The neighbourhood of a cut-element element $\omega_j$ is defined as the set of fluid elements $\{\omega_i\}$ which share a common node with $\omega_j$. This approach is motivated by the extended B-splines proposed in~\cite{Hollig:2001aa}.  It ensures that the polynomial reproduction property of the basis functions is automatically retained so that optimal convergence rates can be achieved irrespective of the extension process. In effect, this technique  enlarges the size of the element from which the basis functions are extrapolated akin to classical element agglomeration strategies. The basis functions with a very small physically active support size and their coefficients do not appear in the system of discretised discontinuous Galerkin equations~\eqref{eq:matrixSys}.   In a cut-element $\omega_j$ with the largest neighbouring fluid element $ \omega_j^{\ast} $, see Figure~\ref{fig:basisExtension},  the interpolation equations are modified as follows 
\begin{equation} \label{eq:intpCut}
\vec{U}_h(\boldsymbol{x},t)\Big|_{\omega_j} = \sum_{i=1}^K  \left [ U^i(t) \varphi_i(\boldsymbol{x}) \right ]_{_{\omega_j^{\ast}}}, \quad 
\vec{S}_h(\boldsymbol{x},t)\Big|_{\omega_j} = \sum_{i=1}^K \left [ S^i(t)  \varphi_i(\boldsymbol{x}) \right ]_{_{\omega_j^{\ast}}}.
\end{equation}
%
In terms of implementation, if the coverage ratio of a cut-element~$\omega_j$ is below $f_j < 0.3$, the discretised conservation vector~$\vec U_h$ and auxiliary vector~$\vec S_h$ on~$\omega_j$ are extrapolated from the basis functions and coefficients of the element $\omega_j^*$. According to~\eqref{eq:intpCut}, this is accomplished by simply evaluating the  interpolation equations for the element~$\omega_j^{\ast}$ at~$\vec x \in \omega_j$. Similarly, the test functions in the weak form~\eqref{eq:weakFormulation} are replaced with the basis functions of the  element $\omega_j^*$. In computing the fluxes~\eqref{eq:LaxFriedrichs},~\eqref{eq:viscousFlux} and~\eqref{eq:auxiliaryFlux} the extrapolated field variables are used.  The explicit time integration of the stabilised system of equations~~\eqref{eq:matrixSys} proceeds as usual. 

In the second proposed stabilisation approach, as introduced in the work of M\"{o}\ss{}ner and Reif~\cite{Mosner2008} for immersed B-spline finite elements, the cut-basis functions are rescaled in proportion to the ratio $1/f_j$.  This leads in a cut-element $\omega_j$  to the following modification of the interpolation equations  
\begin{align}
\vec{U}_h(\boldsymbol{x},t)=\frac{1}{f_j}\sum_{i=1}^K U^i(t)\varphi_i(\boldsymbol{x}), \quad 
\vec{S}_h(\boldsymbol{x},t)=\frac{1}{f_j}\sum_{i=1}^K S^i(t)\varphi_i(\boldsymbol{x}). 
\label{eq:scaledSolutions} 
\end{align}
The implementation of this stabilisation approach is slightly easier than the extrapolation approach. In the field discretisation~\eqref{eq:approximationForm} and the weak form~\eqref{eq:weakFormulation} the basis functions are replaced with the scaled basis functions $\varphi_i(\vec x) / f_j$. This leads to a scaling of the mass matrix~$\ary L$ and the coefficient arrays~$\ary U$ and~$\ary S$ in the semi-discretised  system of equations~$\eqref{eq:matrixSys}$, but not of the right-handside vectors~$\ary R_{\ary U}$ and~$\ary R_{\ary S}$. 
Again, the explicit time integration of the stabilised system of equations~\eqref{eq:matrixSys} proceeds as usual.

\subsection{High-order boundary reconstruction}
%
The cut-element integration using (straight-sided) simplices leads to large errors when the boundary is highly curved. Even when high-order basis functions are used, the domain approximation errors will dominate resulting in a low-order discretisation scheme~\cite{Ciarlet:2002}. In practical terms, an additional complication  concerns the  interpretation of the boundary conditions which are usually prescribed on the curved boundary. 

\begin{figure}[]
\centering
    \subfloat[]{\includegraphics[width=0.45\textwidth]{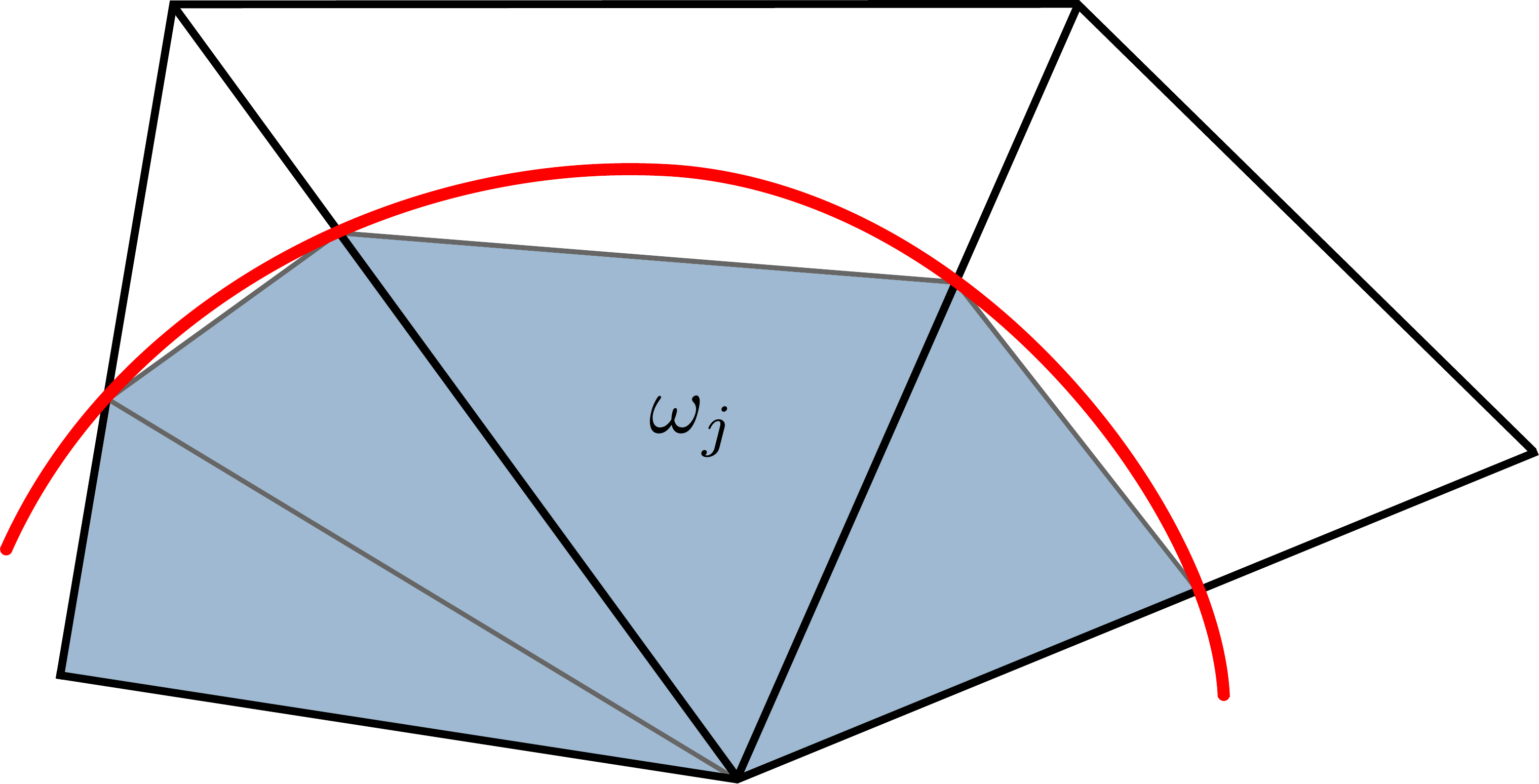}}
    \subfloat[]{\includegraphics[width=0.25\textwidth]{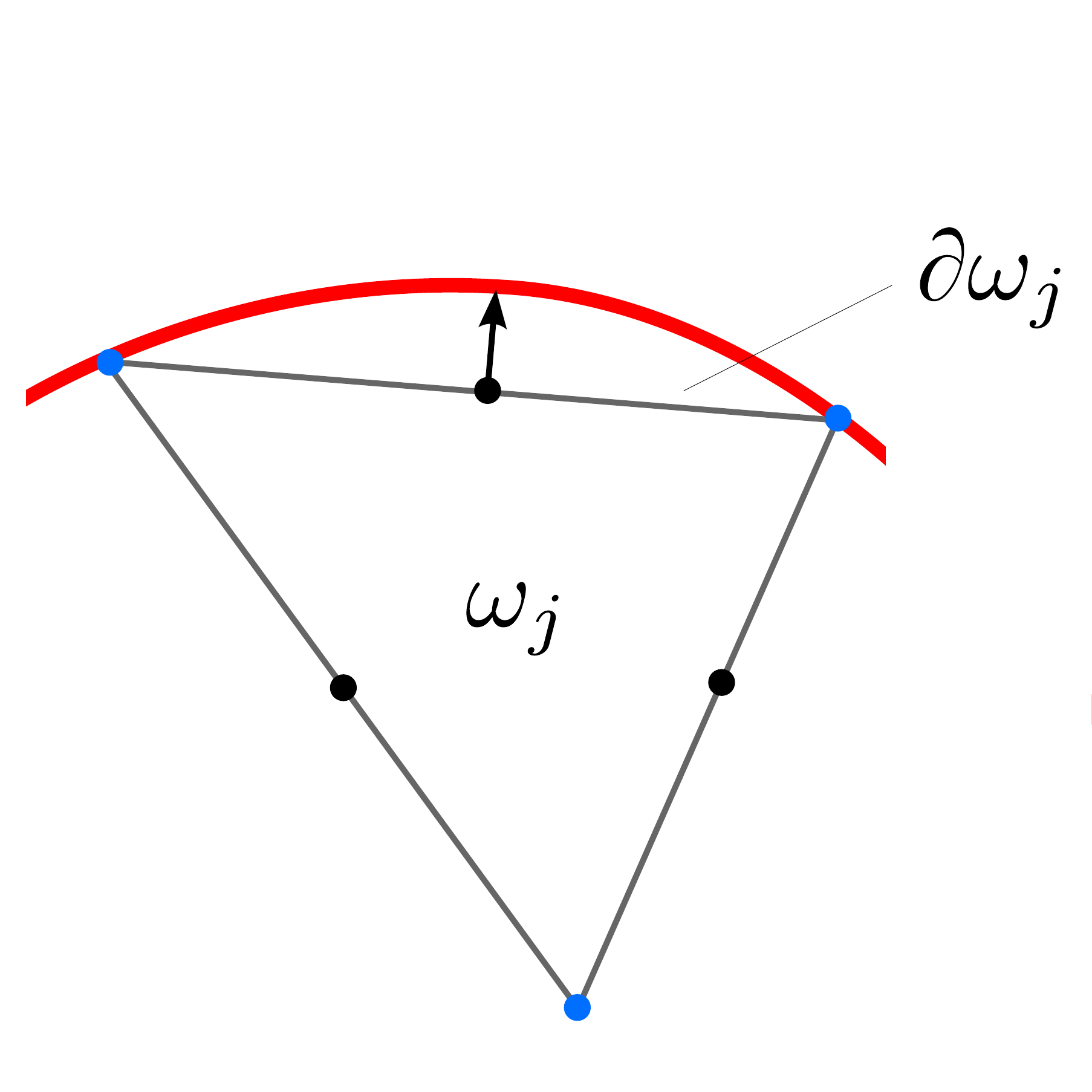}}
    \subfloat[]{\includegraphics[width=0.25\textwidth]{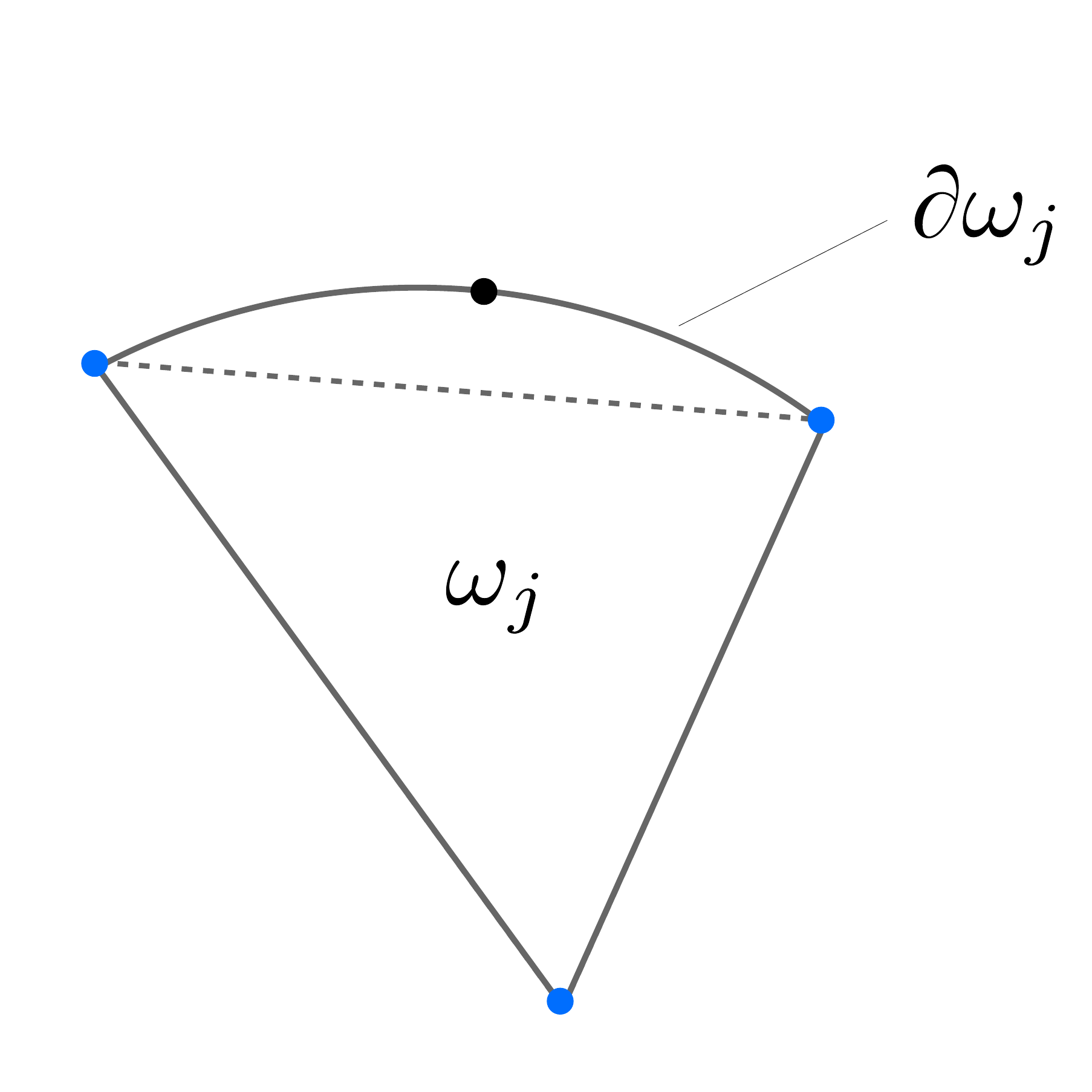}}

  \caption{(a) High-order treatment of a cut-element $\omega_j$ traversed by the curved boundary. (b) Insertion of mid-nodes and their projection to the curved boundary using the level set function~$\phi$ and its gradient $\nabla \phi$. (c) High-order approximation with the domain integrals evaluated on the cut-element with the straight boundary and the boundary integrals evaluated on the curved boundary. }
\label{fig:higherOrderConst}
\end{figure}

It is straightforward to achieve high-order accuracy by suitably curving the boundaries of the simplicial integration elements introduced in Section~\ref{sec:cutElemIntegration} and evaluating all element integrals on the curved simplices.  As shown in the work of Zhang~\cite{Zhang2016}, high-order accuracy can also be achieved by evaluating the domain integrals on straight-sided simplicial elements and evaluating only the boundary, i.e. flux, integrals along curved boundaries, see Figure~\ref{fig:higherOrderConst}.  That is, only the integration of the boundary integrals needs special treatment in case of curved domain boundaries. According to~\cite{Zhang2016}, for nonconvex domains it is assumed that the solution can be smoothly extended outside the domain and that the distance between the straight and the curved boundary is sufficiently small for the sake of stability. Aside from these restrictions, to obtain an optimally convergent discretisation scheme the curved domain boundary  has to be approximated with the same polynomial order like the basis functions used for field interpolation, i.e.~\eqref{eq:approximationForm}. To this end, we modify the boundaries of the integration simplices abutting the domain boundary as follows. Firstly, additional nodes are introduced on the straight-sided simplices. Subsequently, all boundary nodes  are projected to the prescribed domain boundary with the help of the  gradient of the signed distance function~\eqref{eq:sgndistNorm}.  The gradient of the signed distance function is well defined close to a smooth domain boundary. When evaluating the boundary integrals the nonlinear mapping between the straight and curved element boundaries has to be taken into account. To establish this, the curved boundary is assumed to be parametrised with Lagrange polynomials. Throughout the mapping the parametric coordinates remain the same while their physical coordinates change.

As an example, in Figure \ref{fig:higherOrderConst} the process of approximating the domain boundary geometry with quadratic polynomials ($P^2$) is illustrated. First additional mid-nodes are introduced on the edges of the integration simplices. Subsequently, the introduced mid-nodes are pushed along the gradient of the signed distance function~$\nabla \phi$  until the boundary is reached. The intersection of the line given by the gradient with the boundary~$\phi(\vec x)=0$ is found by bisectioning.
%
\section{Results and Discussion}
\label{sec:resultDiscussion}
%
The  theoretical convergence  properties and accuracy of the discontinuous Galerkin method  are well understood for  hyperbolic and elliptic systems. However,  the verification and validation of the discontinuous Galerkin solutions for complex engineering flow problems remain challenging. We present several examples of increasing complexity which are all discretised with unstructured triangular or tetrahedral meshes. As basis functions linear ($P^1$) or quadratic   ($P^2$) Dubiner basis functions are used. The first set of examples verifies the convergence rate of the method and demonstrates the efficacy of the two implemented cut-element stabilisation techniques. In the later more complex examples, our results are both verified and validated mostly by comparison with results from the literature. The last example concerning the hypersonic flow around a vehicle is validated by our experimental results. The examples included cover the subsonic and hypersonic flow regimes. For stationary flows the discontinuous Galerkin equations are advanced in time with explicit time marching until a stationary state is reached. We choose $\CFL=0.8$ for subsonic flows and  $\CFL=0.1$ for hypersonic flows.  The stationary solution is assumed to be reached, depending on the example, either when the increment of the density or drag coefficient  becomes smaller than $10^{-7}$. 

The surface and volume integrals in the weak form  \eqref{eq:weakFormulation} are evaluated numerically. In the present work, the Gauss-Legendre quadrature rule is employed within the elements and element interfaces.  
In two dimensions, we use for the element interfaces 3 Gauss points for linear basis functions ($P^1$),  and 5 Gauss points for quadratic basis functions ($P^2$).  For the element integrals we use 9 and 25 Gauss points  for $P^1$ and $P^2$, respectively.  In three dimensions, boundary integrals are evaluated with 9 and 25 Gauss points for $P^1$ and $P^2$, and volume integrals are  computed with 11 and 21 Gauss points for $P^1$  and $P^2$, respectively.   
\subsection{Flow around a sphere and convergence}
%
We verify the convergence order of the proposed approach with an inviscid flow around a sphere at $Ma=0.2$. The sphere has a radius 1 and is placed at the centre of a cuboidal domain with an edge length of 10. The sphere is defined with the signed distance function 
\begin{equation}
\phi(\boldsymbol{x})=x^2+y^2+z^2-1
\end{equation}
The cuboidal domain is meshed in turn with $64^3$, $96^3$, $ 128^3 $, $ 160^3$ and $ 192^3 $ hexahedrons each of which is later split into six tetrahedra. 

\begin{figure}
	\centering
     \subfloat[]{\includegraphics[width=0.62\textwidth]{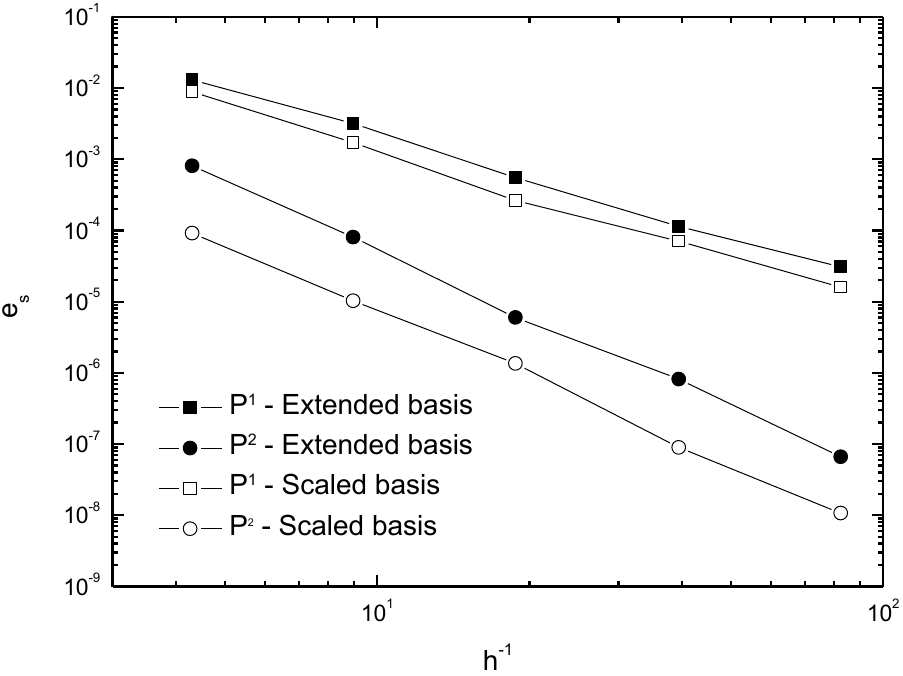}} \hspace{0.025\textwidth}
  	\subfloat[]{\includegraphics[width=0.30\textwidth]{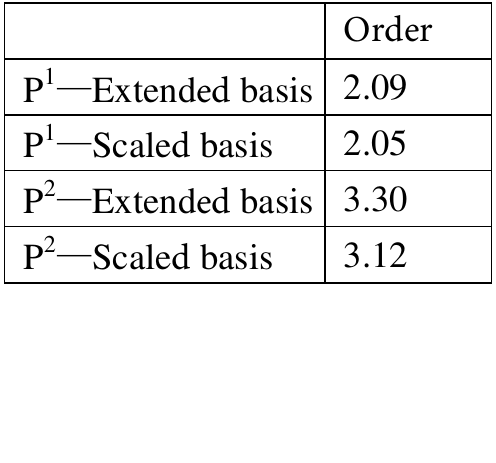}}
    \caption{Flow around a sphere. Convergence of the entropy error  with decreasing characteristic element size~(a). Order of convergence (b). \label{fig:entropyConvergence}}
\end{figure}

The flow is isentropic since it is inviscid and no shocks are present. The numerically computed entropy differs from the free-stream entropy and converges with increasing mesh refinement. The error in the computed entropy  is defined as
\begin{equation}
e_s=\sqrt{\int_\Omega(\Delta s)^2 \, \D v} \, ,
\end{equation}
where $\Delta s=p/\rho^{\gamma}-s_{\infty}$ and the free-stream entropy is $s_{\infty}=p_{\infty}/\rho^{\gamma}_{\infty}=1$.
The convergence of the entropy error with increasing mesh refinement is shown in Figure~\ref{fig:entropyConvergence}. As is to be expected, the entropy error converges with increasing mesh refinement approximately at a rate of three for quadratic ($P^2$) basis functions and  a rate of two for linear ($P^1$) basis functions.  The cut-element stabilisation by extension gives slightly larger errors than by scaling of the basis functions. This may be due to the slight increase of the characteristic element size in the extension approach. 
%

\subsection{Shock trains and cut-element stability}
\begin{figure}
  \centerline{\includegraphics[width=0.9\textwidth]{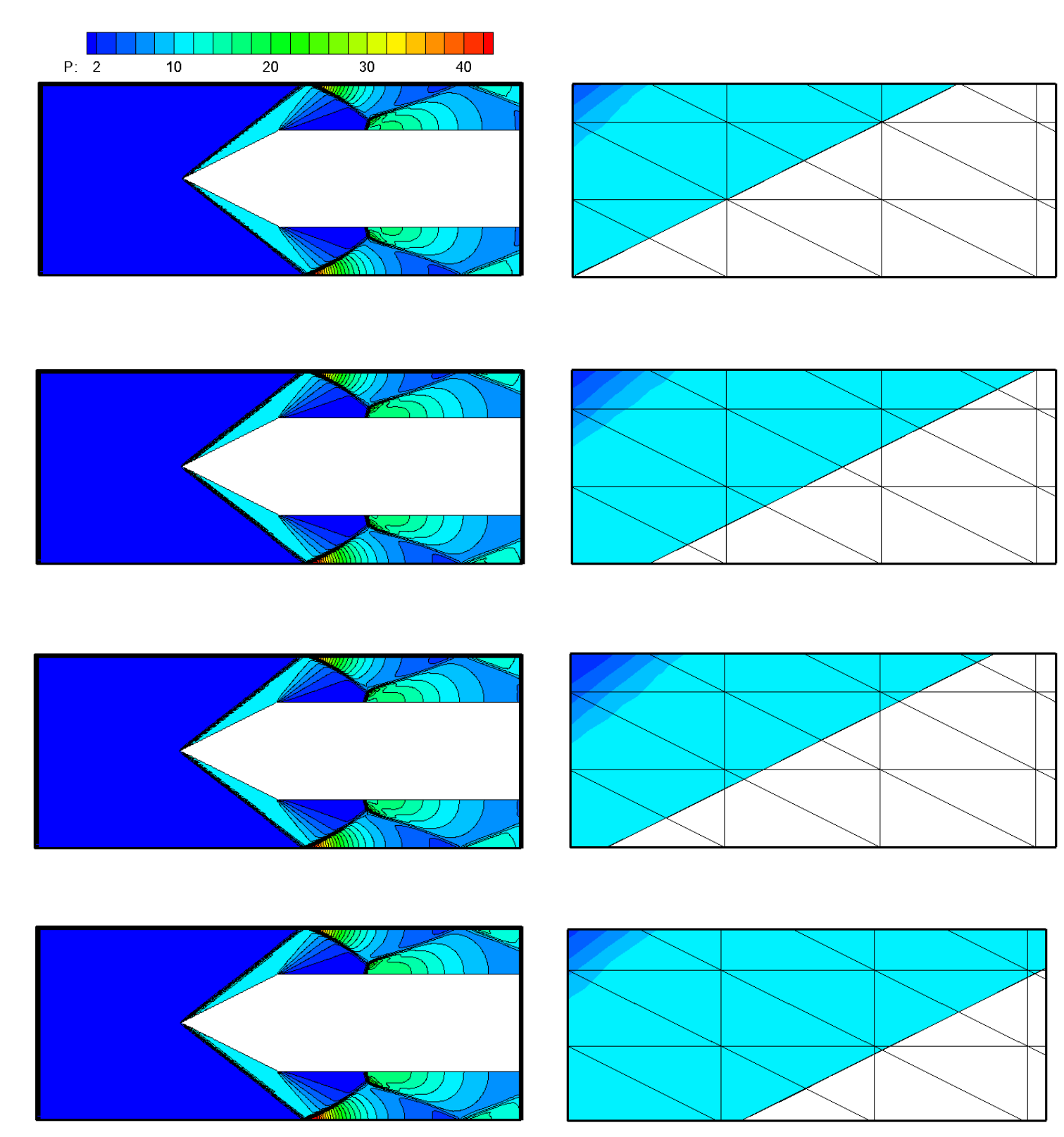}}
  \caption{Shock train pressure isocontours for cut-element stabilisation with the basis scaling approach. Each row shows the pressure isocontour and the close-up of the immersed mesh for the  minimum coverage ratios $f_j \in \{ 1/2, \, 1/8, \, 1/32, \, 1/128 \}$ (from top to bottom).  }
  \label{fig:shockScaling}
\end{figure}

\begin{figure}
  \centerline{\includegraphics[width=0.9\textwidth]{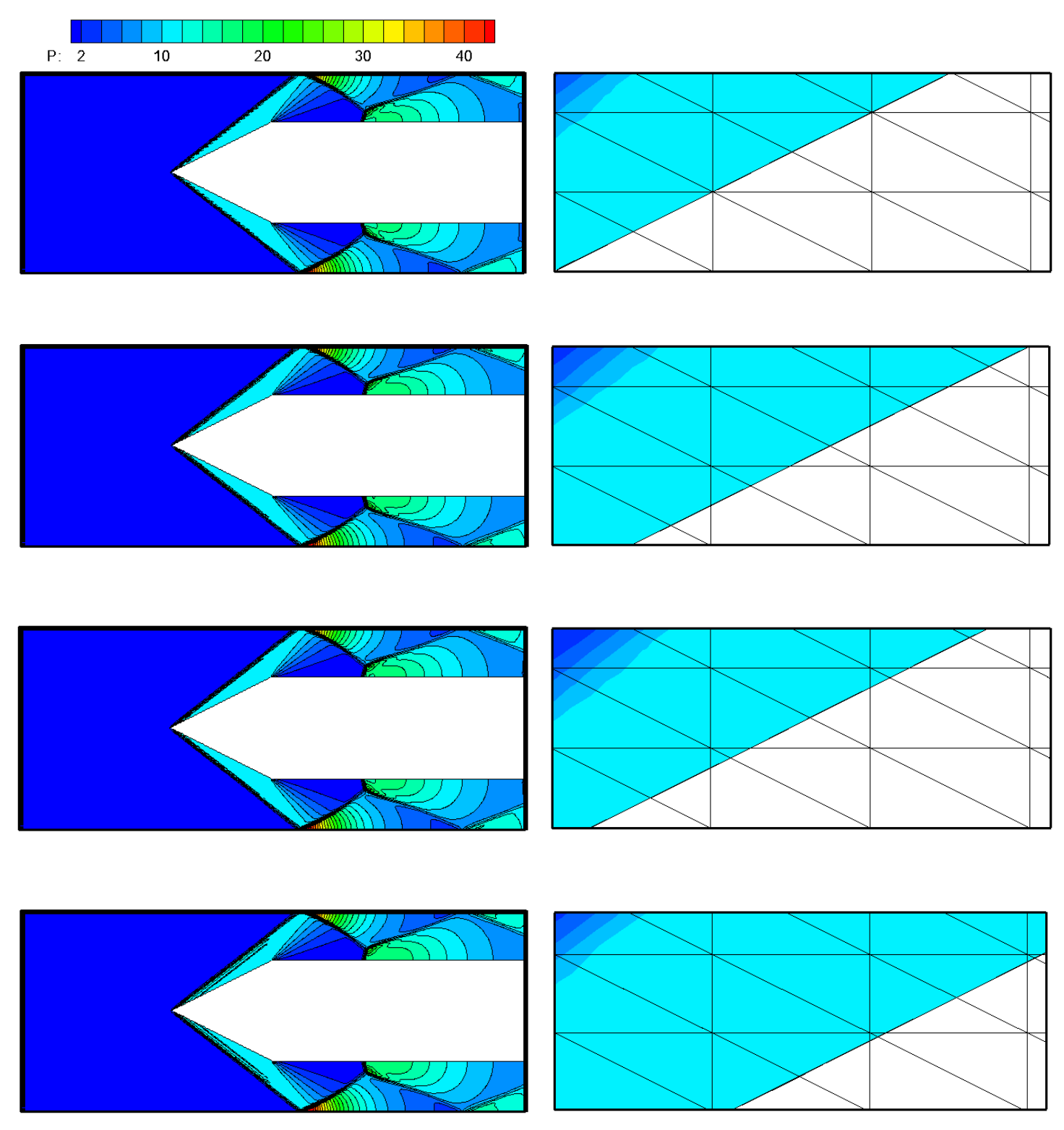}}
  \caption{Shock train pressure isocontours for cut-element stabilisation with the basis extension approach. Each row shows the pressure isocontour and the close-up of the immersed mesh for the  minimum coverage ratios $f_j \in \{ 1/2, \, 1/8, \, 1/32, \, 1/128 \}$ (from top to bottom).}
  \label{fig:shockExtending}
\end{figure}

%
In this example the axisymmetric, two-dimensional compressible Navier-Stokes equations are solved to study the shock train formation inside a confined duct, see Figures~\ref{fig:shockScaling} and~\ref{fig:shockExtending}. The conditions for the flow are $\Mach=5.0$ and $\Reyn=800$.  The Reynolds number is defined with respect to the radius of the solid cone placed inside the duct. The domain is discretised with a structured mesh consisting of  8000 triangular elements  of uniform size and quadratic ($P^2$) Dubiner basis functions.  To study the influence of the proposed cut-element stabilisation techniques the discretisation mesh has been slightly displaced to obtain successively smaller minimum coverage ratios \mbox{$f_j \in \{ 1/2, \, 1/8, \, 1/32, \, 1/128 \}$.}

Stationary shock trains consisting of a series of successive shock and rarefaction waves  are frequently used for verification and validation of high-order schemes ~\cite{Matsuo1999}. Although such problems have been extensively studied, see e.g. \cite{Carroll1990, rodriguez2001asymmetry}, they provide valuable information about the  capability of the proposed immersed discontinuous Galerkin method in dealing with intricate gas dynamics problems. The structure of the shock train can be complex even for simple inflow and wall conditions and the flow remains hypersonic throughout the duct. For instance, in a symmetric duct with uniform, steady inflow conditions, asymmetric flow structures may develop and self-excited oscillations may appear~\cite{rodriguez2001asymmetry,ikui1974oscillation}.  These features are not present in our simplified axisymmetric, two-dimensional set-up. The shock trains are usually classified as either normal or oblique shock trains depending on the presence or absence of bifurcated normal shock waves~\cite{Carroll1990}. Lower upstream Mach numbers typically lead to normal shock trains and higher Mach numbers, as in our computations, lead to oblique shock waves. 

The pressure isocontours of the present computations are shown in Figures~\ref{fig:shockScaling}  and  \ref{fig:shockExtending} for the scaled and extended basis stabilisation techniques, respectively.  As the gas flow impinges on the cone, first a planar oblique shock is generated. After this shock wave impinges on the duct wall, a reflected shock wave is generated. This process continues and a series of reflected shock waves are created.  The results for four sets of meshes and two cut-element stabilisation  techniques show that the size, structure and locations of the incident and reflected shock waves are visually identical.

\subsection{Flow around a two-dimensional cylinder}
\begin{figure}[h]
 \centering
 \subfloat[\label{fig:meshCylinderA}]{\includegraphics[width=0.32\textwidth]{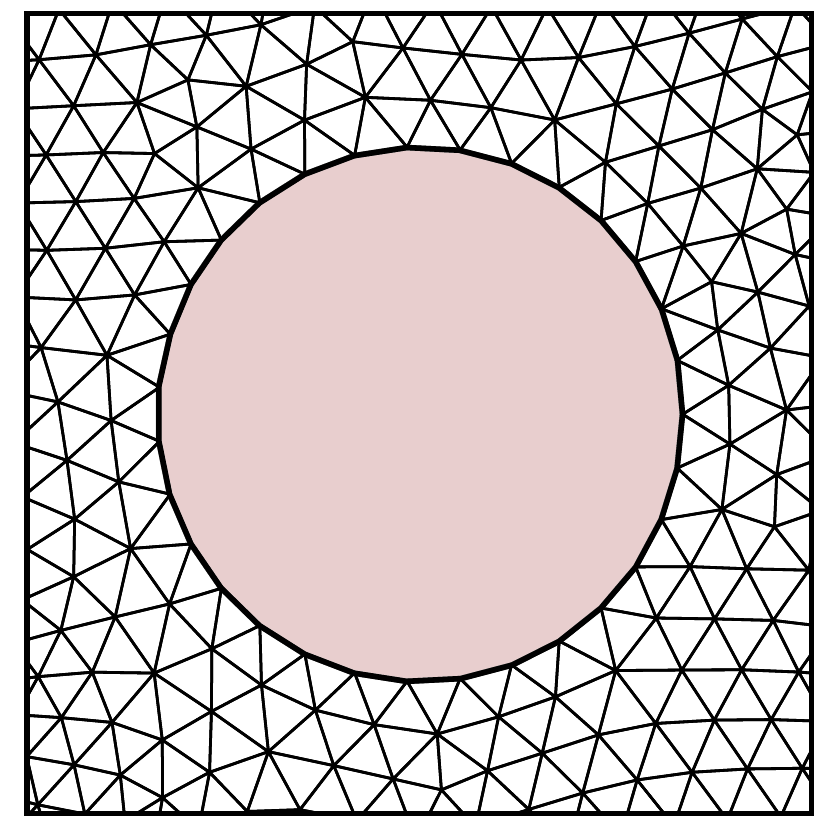}} \hfill
 \subfloat[\label{fig:meshCylinderB}]{\includegraphics[width=0.325\textwidth]{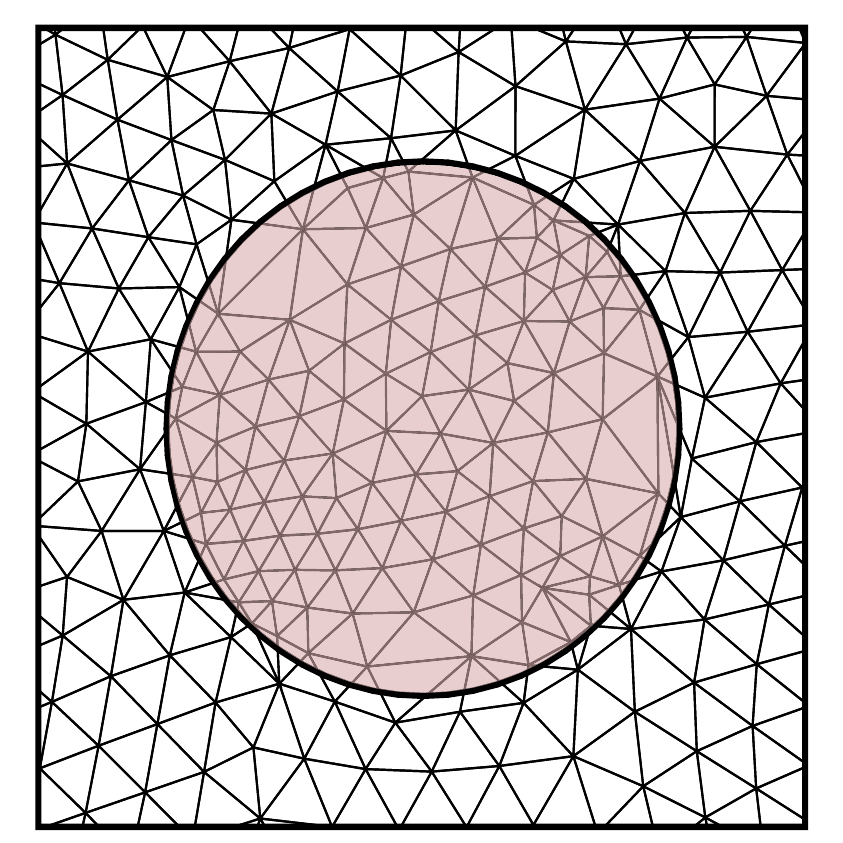}} \hfill
  \subfloat[\label{fig:meshCylinderC}]{\includegraphics[width=0.32\textwidth]{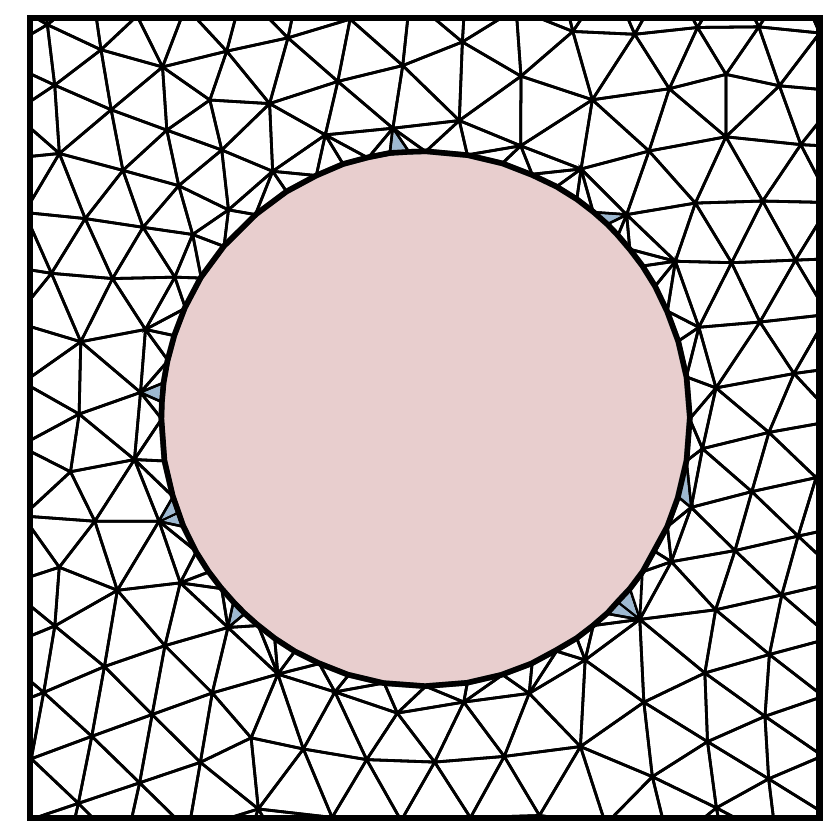}}
     \caption{Close-ups of the boundary-fitted (a) and non-boundary-fitted (b) meshes for the flow around a two-dimensional cylinder. The non-boundary-fitted mesh contains~$62$ cut-elements which are for integration  triangulated with the marching triangle algorithm (c). The nine blue shaded cut-elements require stabilisation.}
  \label{fig:geomCylinder}
\end{figure}

%
There is extensive literature on the study of sub- and supersonic flows around a cylinder, which makes it an ideal application to verify and validate the developed immersed discontinuous Galerkin method. In all the presented computations, a cylinder with diameter $d=1$ is placed at the centre of a fluid domain with radius $30 d$. For meshing, $100$ nodes are placed along the circumference of the outer fluid domain boundary and across the diameter, on average, $300$ nodes are present. In the computations the boundary-fitted and non-boundary-fitted meshes depicted in Figures~\ref{fig:meshCylinderA} and~\ref{fig:meshCylinderB} are used.  Each of the meshes consists of about $60000$ triangular elements. As basis functions the quadratic ($P^2$) Dubiner basis functions are used.  Along the circumference of the outer fluid domain boundary far-field boundary conditions are prescribed. 

\subsubsection{Subsonic flow}
%
For very low Reynolds numbers, with $\Reyn \ll1$ (based on cylinder's diameter), the flow is symmetrical upstream and downstream of the cylinder. As Reynolds number increases, the upstream-downstream symmetry breaks. When the Reynolds number  exceeds $\Reyn > 4$ two attached eddies emerge behind the cylinder, see e.g.~\cite{wieselsberger84new, TrittonD.J.1959, panton2013incompressible}. These eddies increase in size with increasing Reynolds numbers before the flow becomes nonstationary at about $\Reyn \gtrsim 46$ and periodic vortex shedding begins~\cite{Ye1999}. 

\begin{figure}
  \centering 
  \includegraphics[width=0.7\textwidth]{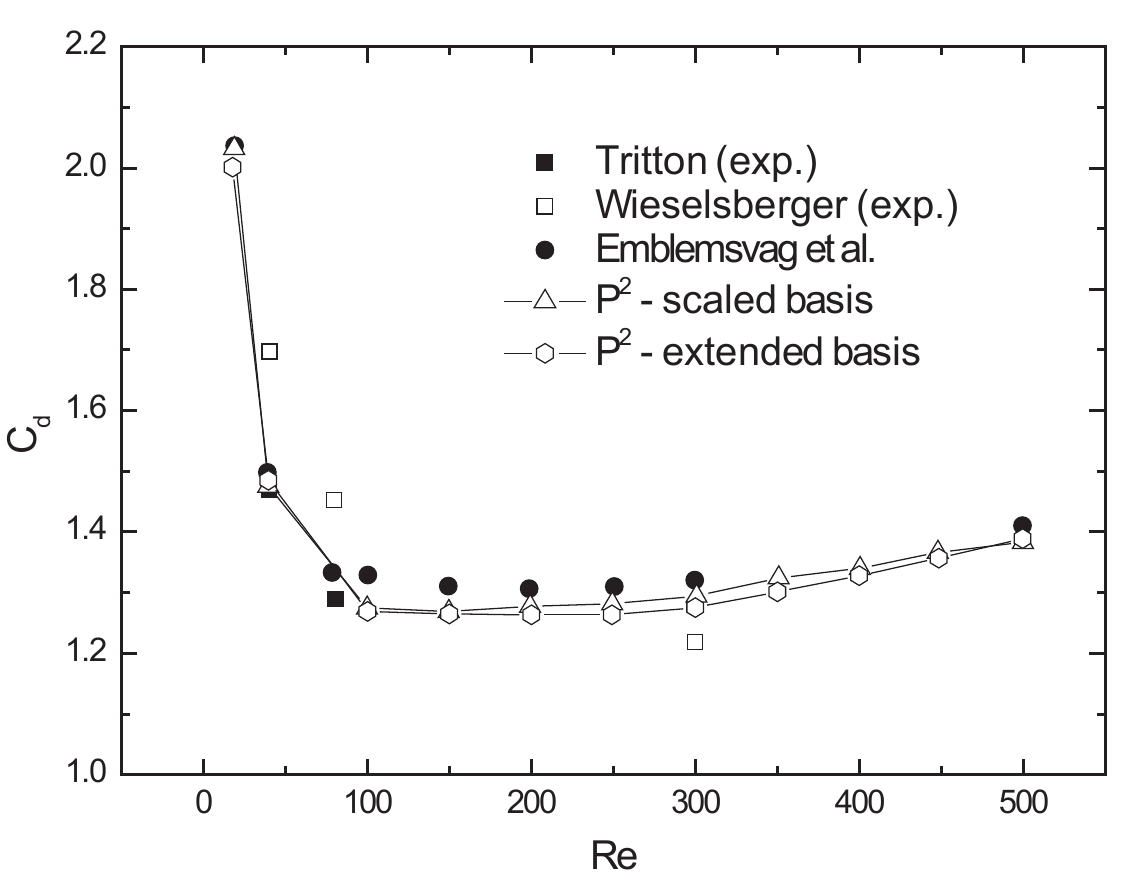}
  \caption{Average drag coefficient $C_d$ of a cylinder with increasing Reynolds number.}
  \label{fig:dragCylinder}
\end{figure}

In the first set of computations, we use the non-boundary-fitted mesh depicted in Figure~\ref{fig:meshCylinderB} to verify and validate the accuracy of the developed techniques. This mesh contains~$62$ cut-elements and nine of them require stabilisation. In Figure~\ref{fig:dragCylinder} the computed average drag coefficients $C_d$ for Reynolds numbers up to $\Reyn = 500$ are plotted. The obtained values are in close agreement with the computational results reported by Emblemsvag et al. \cite{Emblemsvag2005} and the experimental results reported by  Wieselsberger \cite{wieselsberger84new} and Tritton \cite{TrittonD.J.1959}. 
In addition to the average drag coefficient, the wake separation length (for $\Reyn \lesssim 40$)  is also a feature of the flow that is often compared. The wake separation length is defined as the distance from the rear of the cylinder to the tail end of the wake. Table~\ref{tab:wakeSepLen} shows that our results are in very good agreement with computations by others. Moreover, the introduced two cut-element stabilisation techniques lead to almost the same results.

\begin{table}
\caption{Comparison of wake separation length $L_w/d$. }
\centering
\begin{tabular}{lll}
\hline
  $\Reyn$ & 20 & 40\\
\hline
  Dennis and Chang \cite{Dennis2006} & 0.94 & 2.35\\
  Fornberg \cite{Fornberg1980} & 0.91 & 2.24\\
  Ye et al. \cite{Ye1999} & 0.92 & 2.27\\
  Emblemsvag et al. \cite{Emblemsvag2005} & 0.93 & 2.28\\
  $P^2$ - scaled basis & 0.920 & 2.280\\
  $P^2$ -  extended basis & 0.929 & 2.284\\
\hline
\end{tabular}
\label{tab:wakeSepLen}
\end{table}

\begin{figure}
\centering 
 \includegraphics[width=1\textwidth]{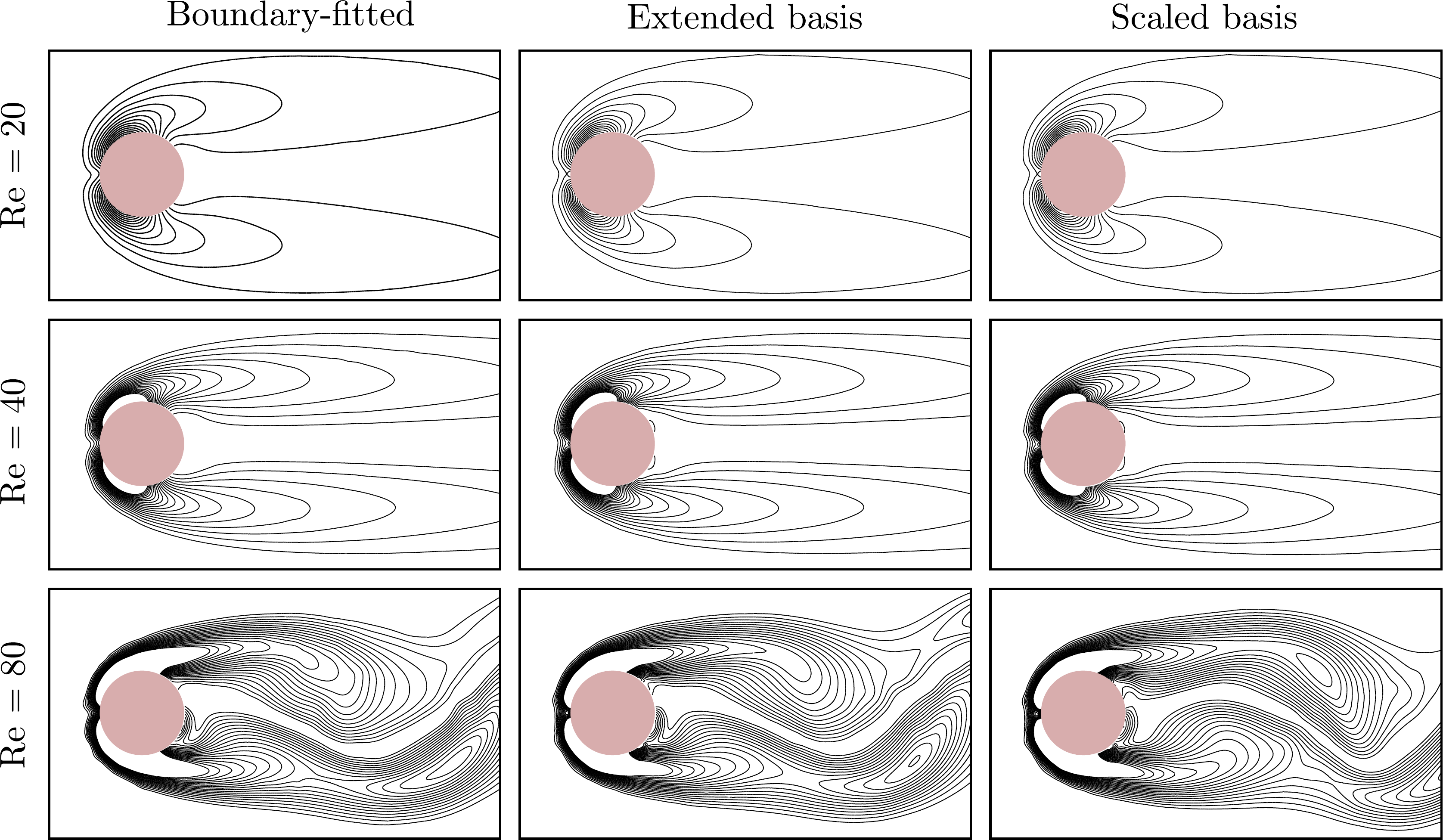}
 \caption{Vorticity isocontours for the flow around a cylinder at~$\Reyn=20$, $\Reyn=40$  and~$\Reyn=80$ (from top to bottom). The isocontours range from~$-4$ to~$4$ with an interval size of~$0.2$.  The left column gives the results for the boundary-fitted mesh and the other two columns give the results for the non-boundary-fitted mesh stabilised either by basis extension (middle) or scaling (right). 
  \label{fig:vorticityIsocontour} }
\end{figure}

In the second set of computations, we compare the vorticity fields obtained with the boundary-fitted and non-boundary fitted meshes. The vorticity field is highly sensitive to the accurate resolution of the velocity gradients in the vicinity of the cylinder boundary. In Figure~\ref{fig:vorticityIsocontour} the vorticity isocontours for the boundary-fitted and non-boundary fitted meshes and three different Reynolds numbers~$\Reyn=20$, $\Reyn=40$ and~$\Reyn=80$ are plotted. 
For the two lower Reynolds numbers~$\Reyn=20$ and~$\Reyn=40$ the flow is stationary. As evident from the vorticity isocontours, the boundary-fitted and non-boundary fitted meshes yield very similar results. For~$\Reyn=20$ the three isocontour plots are essentially indistinguishable. For~$\Reyn=40$ there is some difference in the back of the cylinder, which appears to be due to the slightly larger recirculation velocity in case of the immersed method. For the instationary~$\Reyn=80$ the essential features of the vorticity field are the same while there are some minor differences in the wake flow. In making these comparisons it has to be noted that the element sizes and connectivity of the boundary-fitted and non-boundary-fitted meshes are, especially close to the cylinder, different. 

\subsubsection{Hypersonic flow}
%
Next, the hypersonic flow around  a two-dimensional cylinder is considered  to validate the proposed immersed discontinuous Galerkin method.  Shock detection and capturing in high-order discontinuous Galerkin methods is challenging and remains an active area of research.  For finite difference and finite volume methods many ingenious shock capturing methods are known. Although the extension of these methods to the one-dimensional  discontinuous Galerkin method works well, their extension to two- and three-dimensional unstructured meshes is not always satisfactory. 

We consider the two experiments reported by Maslach et al.~\cite{Maslach1962}. Between the two experiments only the Knudsen number $\Knud=0.016$ and $\Knud=0.08$  is different. The other flow parameters are $\Mach=5.92$,  $\Pran = 0.707$, $T_{\infty}=36.6K$ and $T_w=293K$. The dimensionless parameters correspond to a physical problem with cylinder radius $R=1.55\, mm$ and $R = 0.31\, mm$ with air as the working gas.

\begin{figure}
\centering
    \subfloat[$\Mach=5.92, \Knud = 0.016$]{\includegraphics[width=0.4\textwidth]{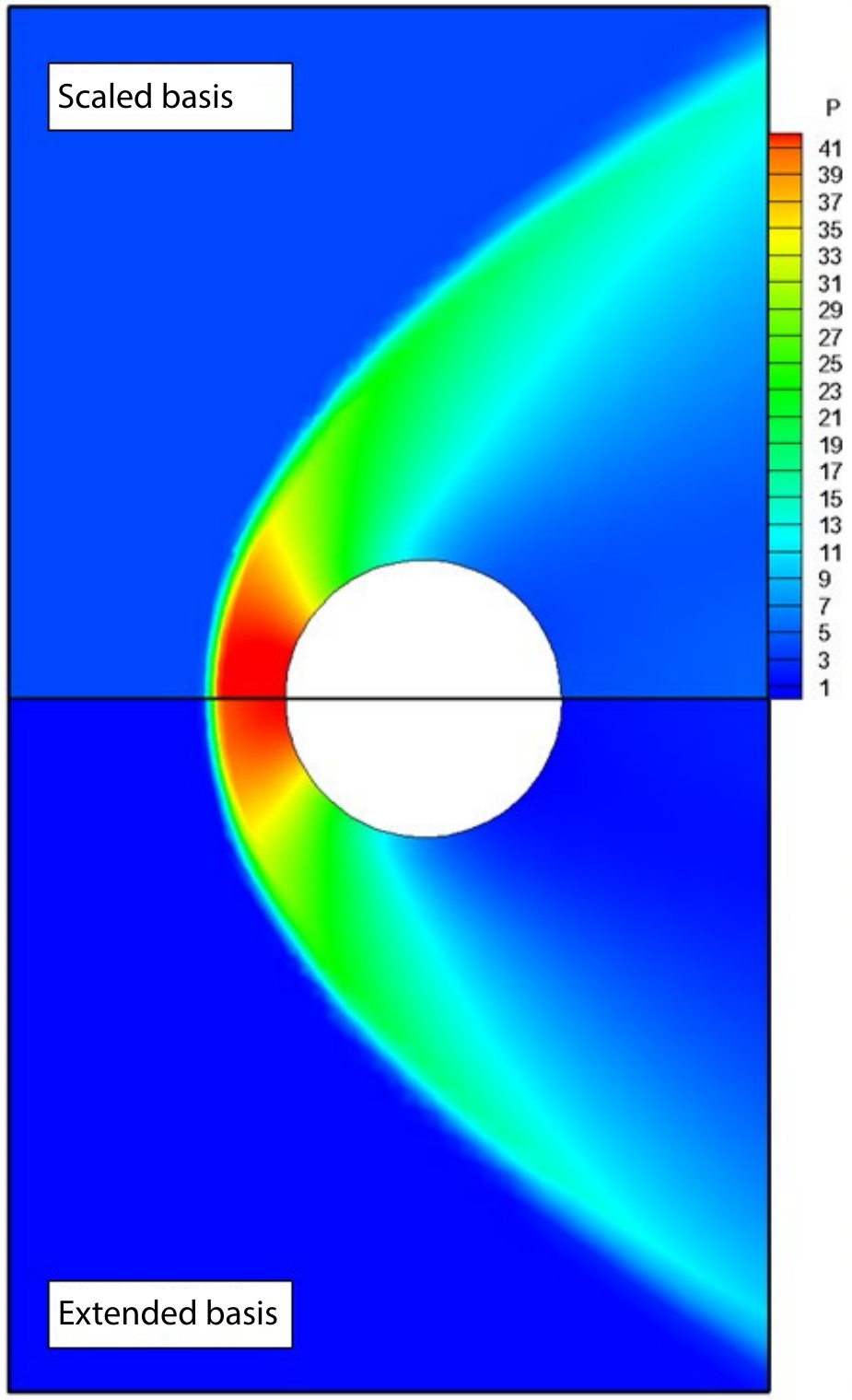}}
    \hspace{0.1\textwidth}
	\subfloat[$\Mach=5.92, \Knud = 0.08$]{\includegraphics[width=0.4\textwidth]{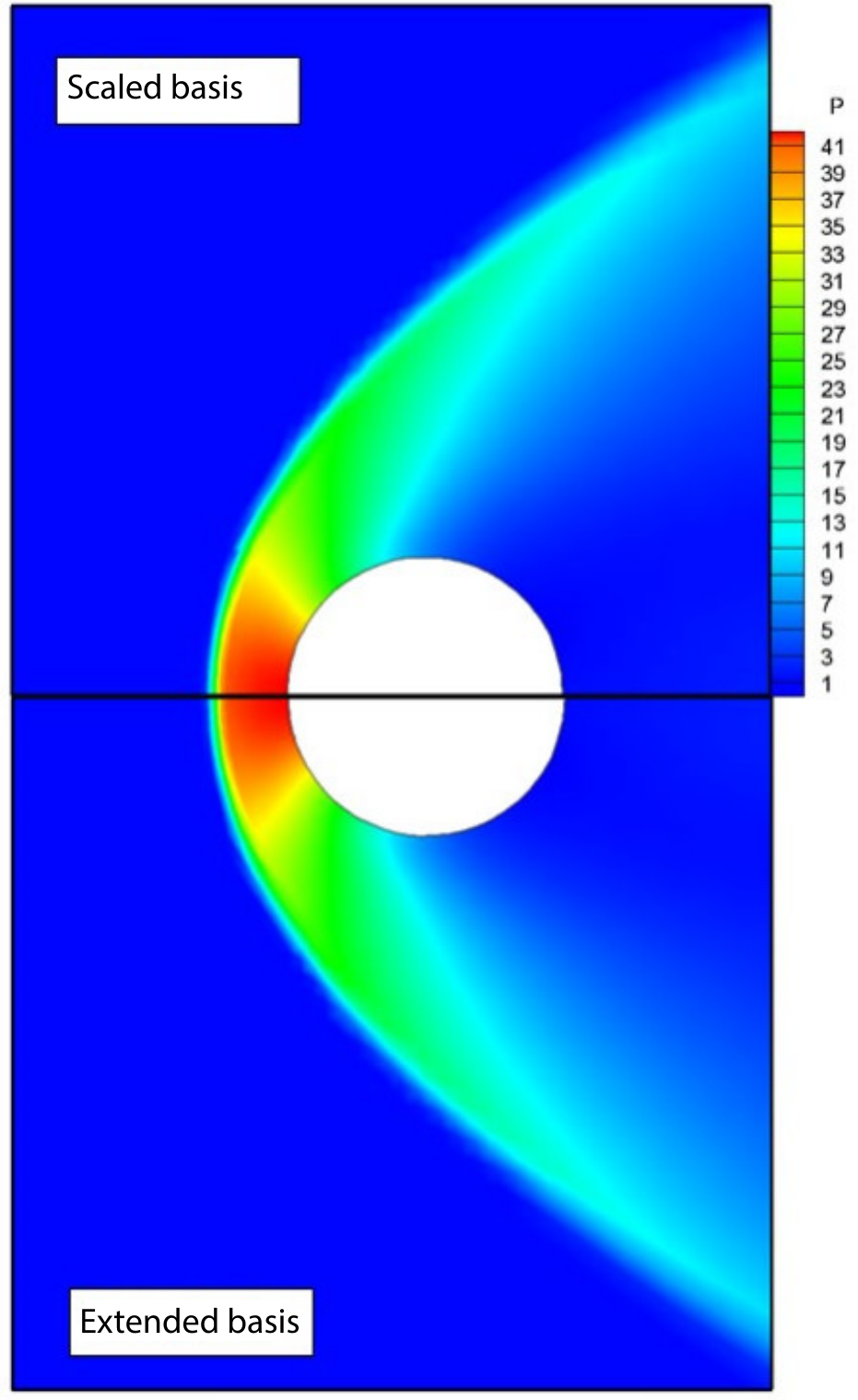}}
  \caption{Pressure isocontours of the hypersonic case.}
  \label{fig:hyperCylinder}
\end{figure}

In Figure~\ref{fig:hyperCylinder} the computed pressure isocontours for the two experiments are shown. Both computations have been in turn computed with the scaled and extended basis stabilisation techniques.  As evident from  Figure~\ref{fig:hyperCylinder} the two alternative cut-element stabilisation techniques yield indistinguishable results in terms of bow shock standoff distance and pressure contours in the proximity of the cylinder. In Table~\ref{tab:dragCompare} the computed drag coefficients are reported which are in excellent agreement with the experiments, with the scaled basis stabilisation giving slightly better results.

\begin{table}
\caption{Comparison of the computed drag coefficients $C_d$ with the experiments.}
\centering

\begin{tabular}{lll}
\hline
  \quad & $\Mach$=5.92, $\Knud$=0.016 & $\Mach$=5.92, $\Knud$=0.08 \\
\hline
  Exp. by Maslach et al. \cite{Maslach1962} & 1.55 & 1.63\\
 $P^2$ - scaled basis (error) & 1.561 (0.71\%)  & 1.647 (1.04\%) \\
  $P^2$ - extended basis  (error) & 1.563 (0.84\%) & 1.649 (1.17\%) \\
\hline
\end{tabular}
\label{tab:dragCompare}
\end{table}

\subsection{Hypersonic gas flow around the Apollo 6 command module}
\begin{figure}
  \centering
  \includegraphics[width=0.75\textwidth]{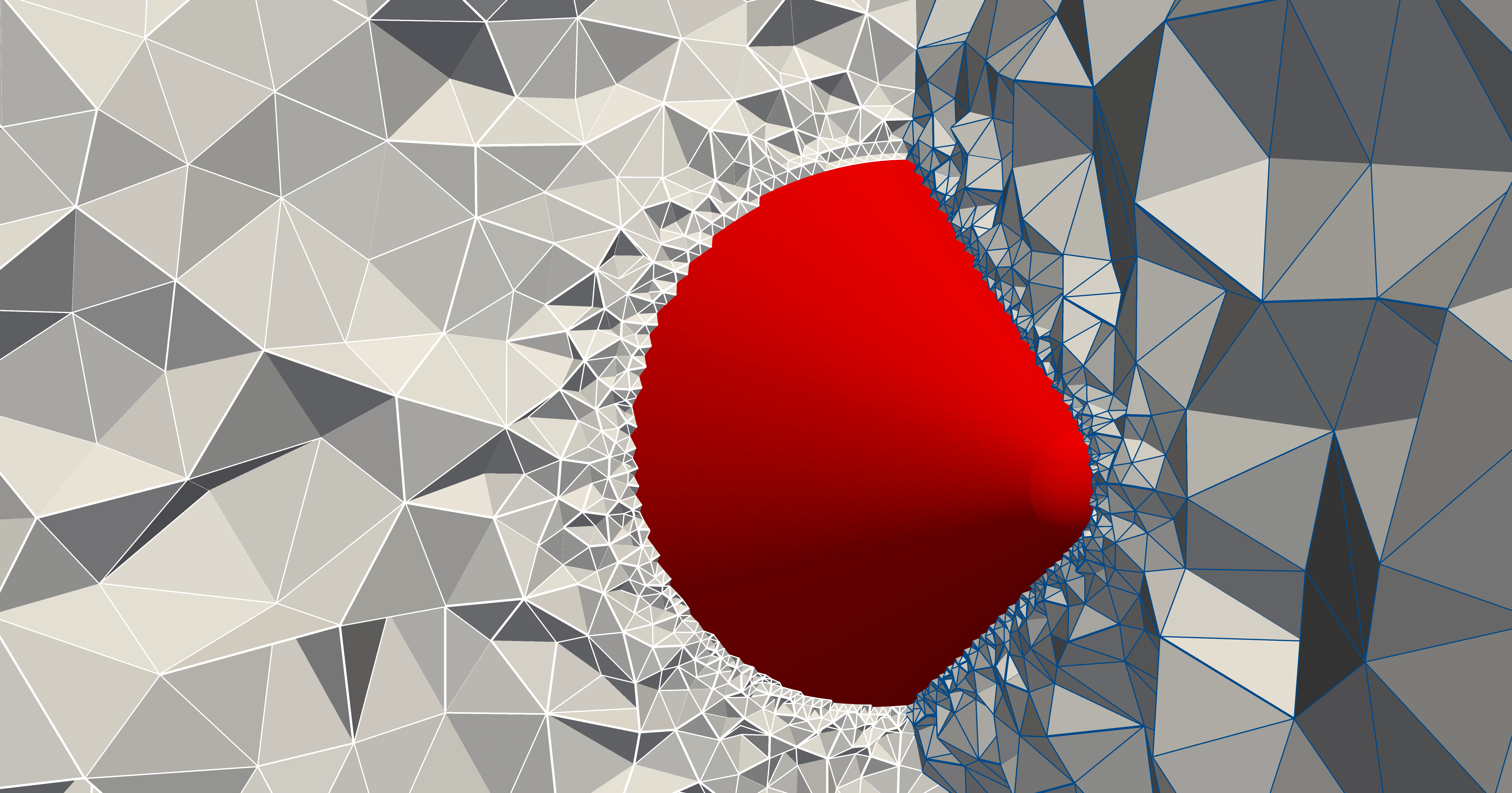}
  \caption{Apollo 6 Command Module geometry and mesh.}
  \label{fig:ApolloGeom}
\end{figure}

\begin{figure}
  \centering 
  \subfloat[Scaled basis]{\includegraphics[width=0.45\textwidth]{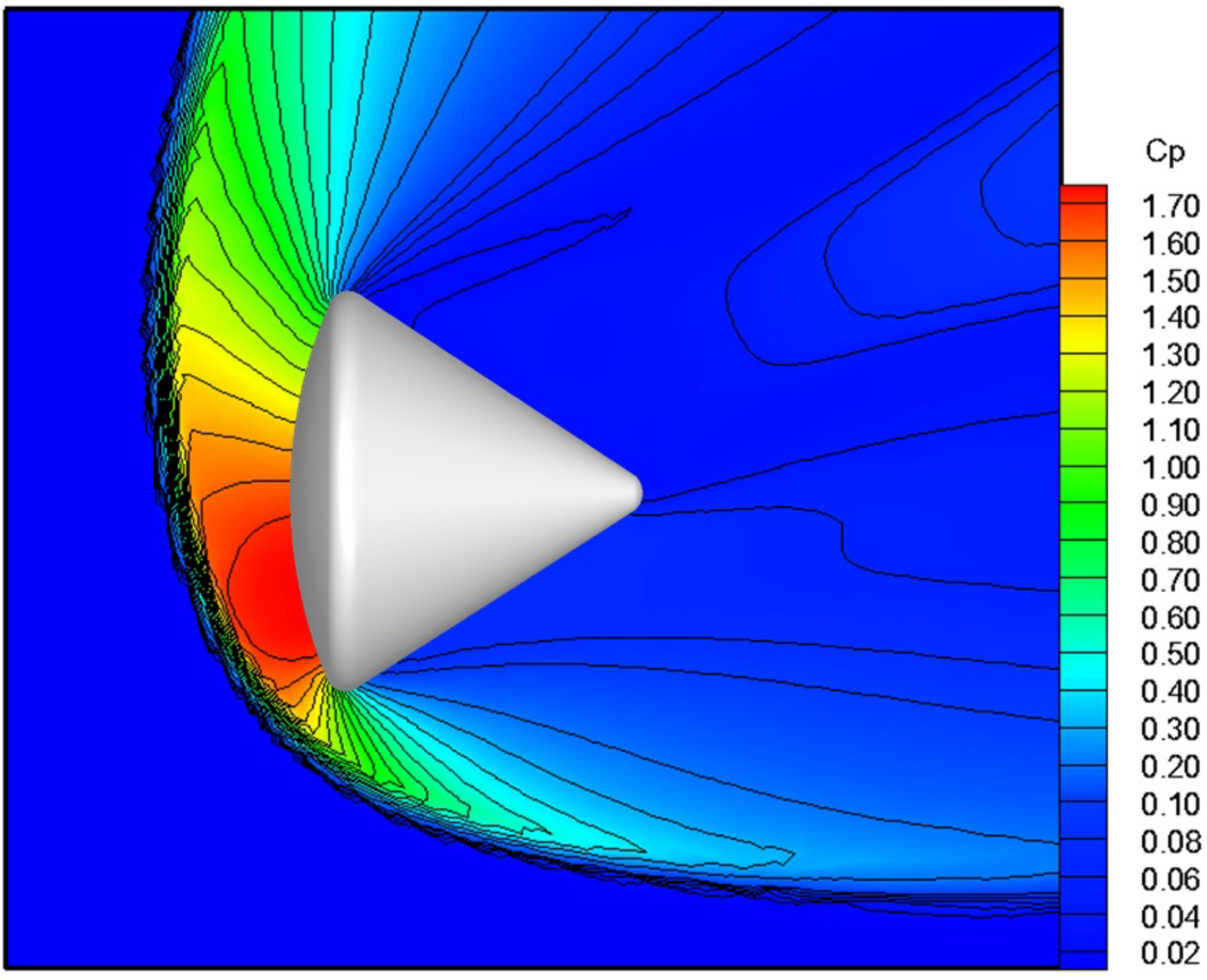}} \hspace{0.05\textwidth}
   \subfloat[Extended basis]{\includegraphics[width=0.45\textwidth]{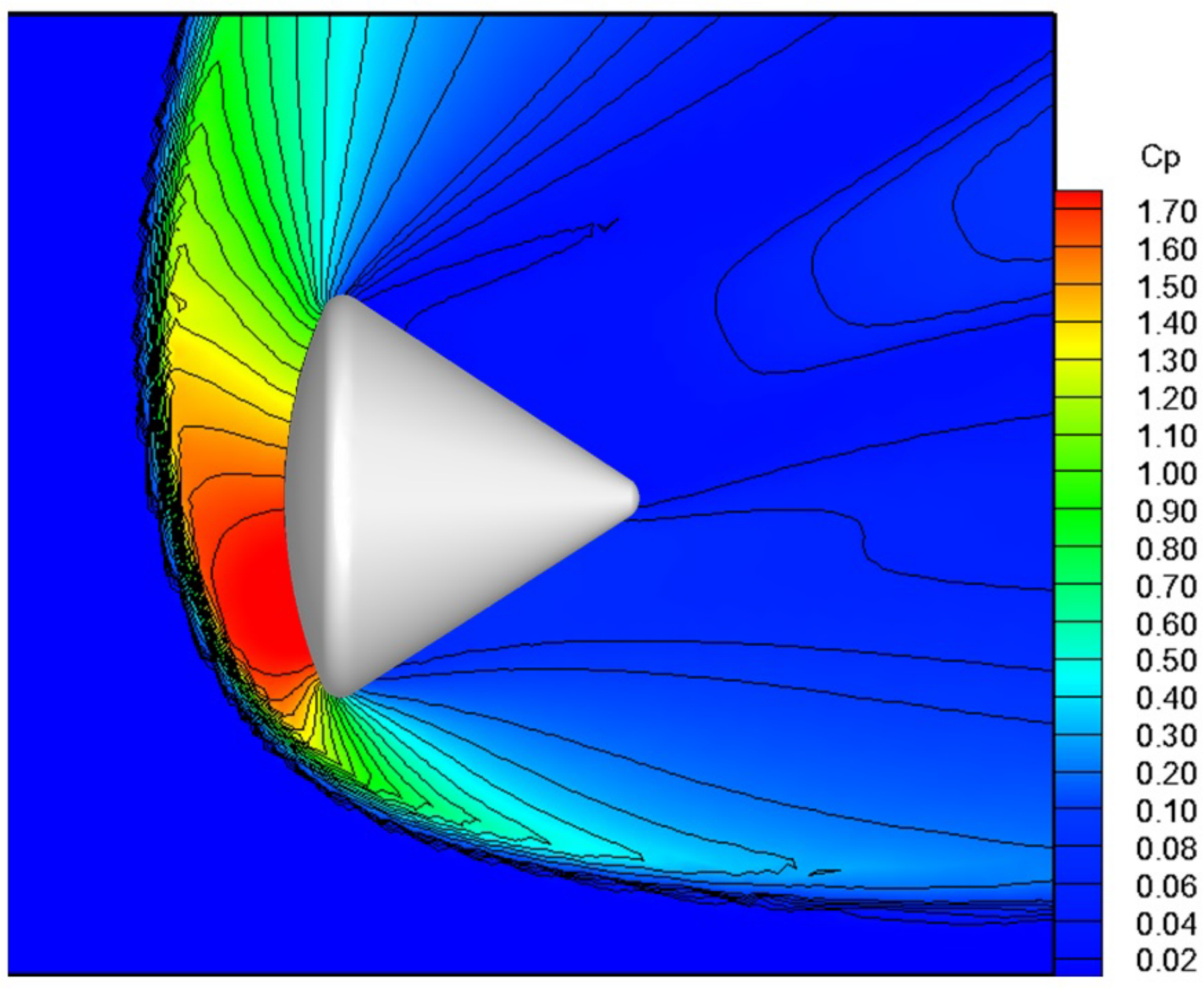}}
  \caption{Pressure isocontours of the flow around the Apollo 6 Command Module.}
  \label{fig:ApolloRes}
\end{figure}

%
As a more complex engineering problem, the three-dimensional gas flow around the Apollo 6 command module \cite{Moss2007} is investigated. The assumed freestream flow and wall conditions are $\Mach=5.0$,  $\Reyn = 1500$, $T_{0}=288 \,K$ and $T_w= 288 \, K$.  The angle of attack is $25^\circ$. The axisymmetric command module geometry and the unstructured tetrahedral discretisation mesh are shown in Figure~\ref{fig:ApolloGeom}. The command module is placed in the centre of a spherical fluid domain with a diameter of $30 \,d$, where $d$ is the maximum diameter of the module. The mesh consists of 386600 elements in total and across the diameter of  the fluid domain, on average, 720 elements are present. As basis functions the quadratic ($P^2$) Dubiner basis functions are used and the total number of degrees of freedom is approximately 19.3 million. The surface mesh generated during cut-element integration consists of 54576 triangular elements.  The computed pressure isocontours are shown in Figure~\ref{fig:ApolloRes} with the  clearly captured shockwave. Unfortunately, we do not consider the reacting flow so that it is not possible to compare our computational results with experimental data. 

\subsection{Hypersonic gas flow around a complex vehicle}
\begin{figure}[]
  \centerline{\includegraphics[width=0.8\textwidth]{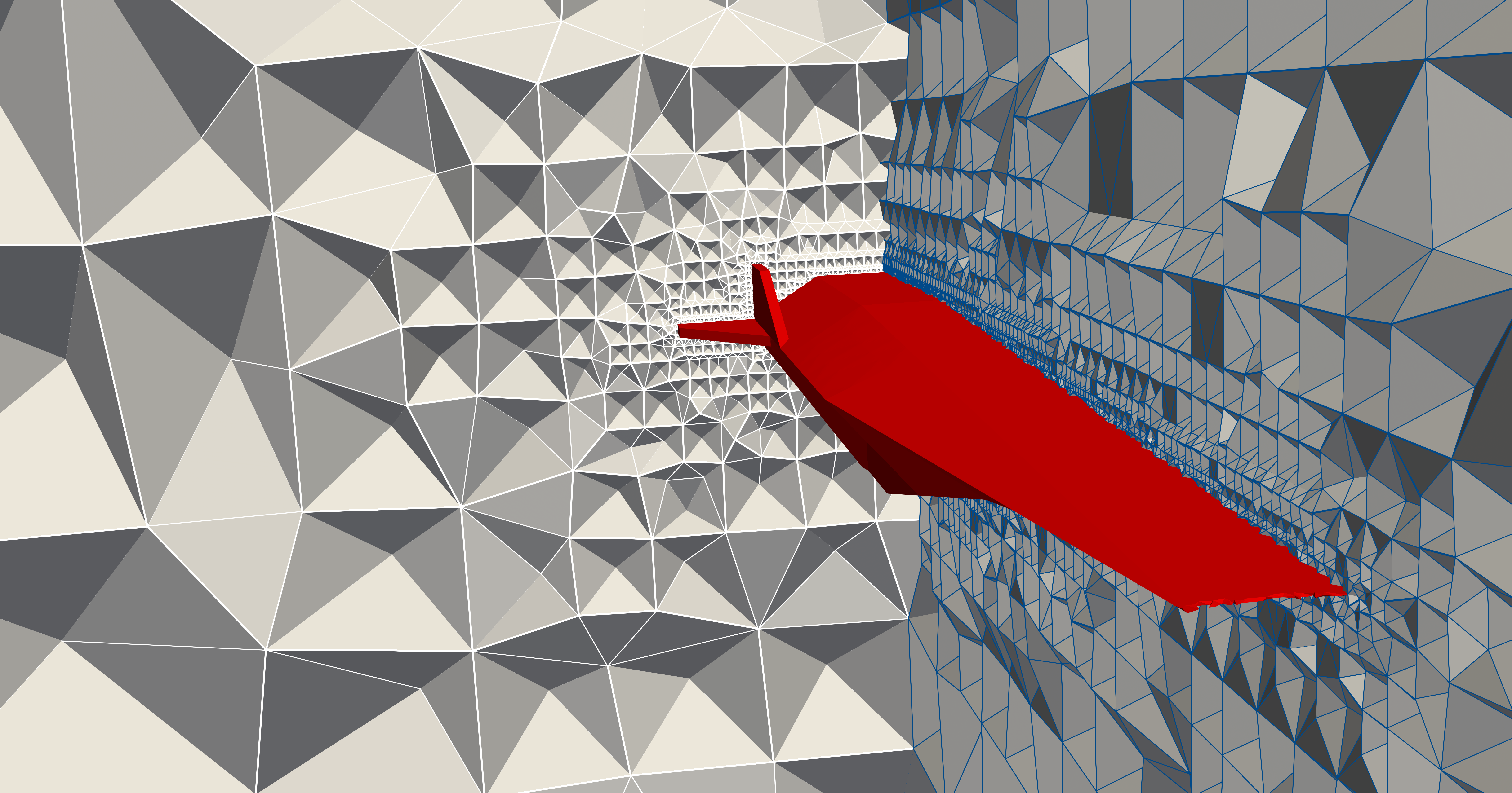}}
  \caption{Hypersonic air vehicle geometry and mesh. }
  \label{fig:geomHypersonic}
\end{figure}

\begin{figure}[]
  \centerline{\includegraphics[width=0.75\textwidth]{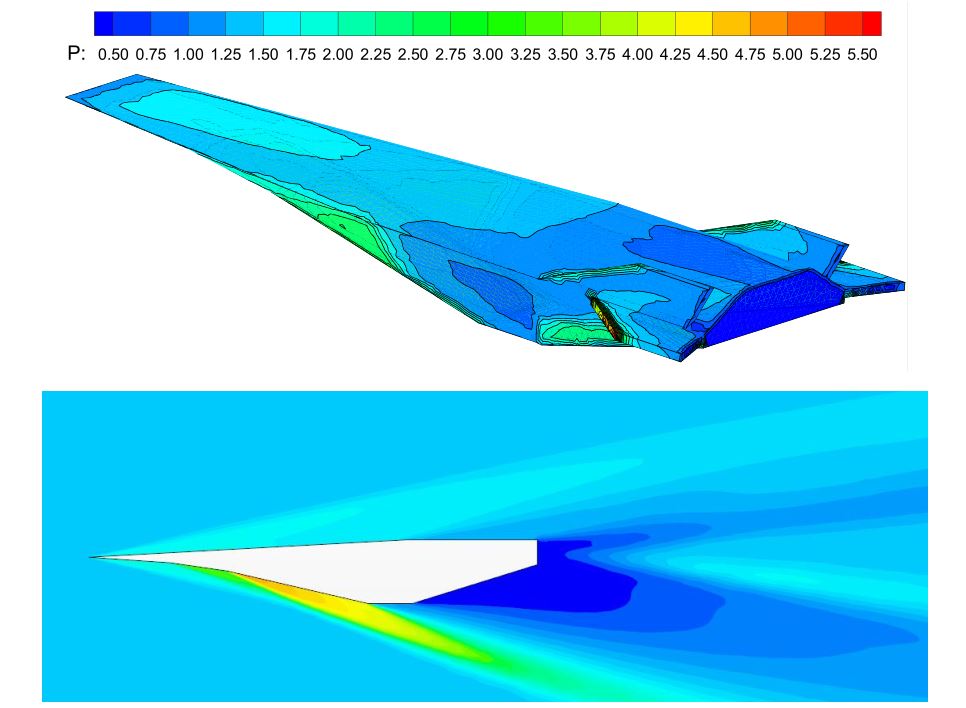}}
  \caption{Pressure iso-contours on the surface of the hypersonic air vehicle (top) and a cross section of the flow domain (bottom).}
  \label{fig:pressHypersonic}
\end{figure}

\begin{figure}
  \centerline{\includegraphics[width=0.7\textwidth]{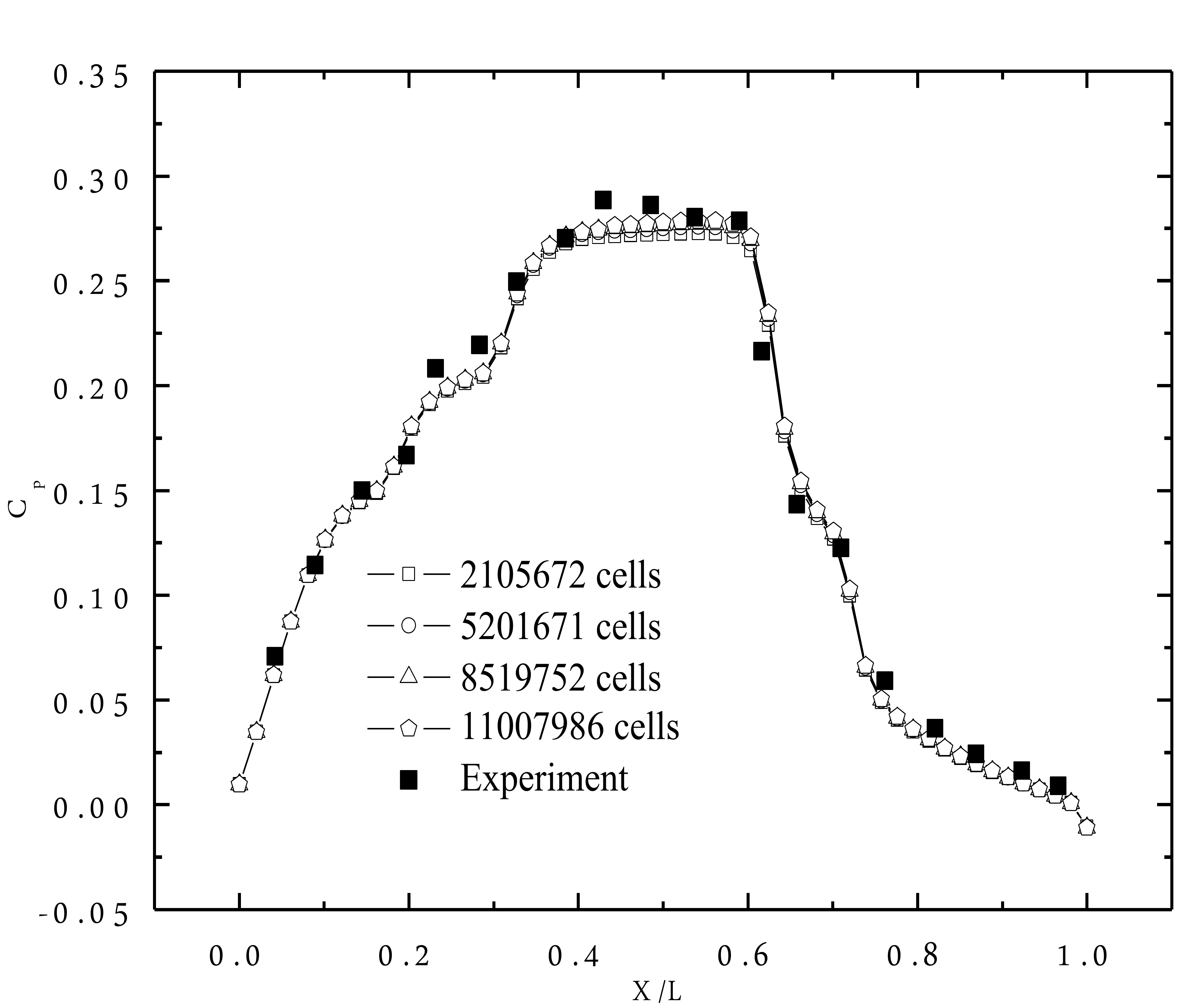}}
  \caption{Pressure coefficient $C_p$ along the symmetry axis on the bottom surface of the hypersonic air vehicle.}
  \label{fig:pressCoefHypersonic}
\end{figure}

%
As the final example, the complex three-dimensional gas flow around a hypersonic aerospace vehicle is considered, see Figure~\ref{fig:geomHypersonic}.  The freestream flow and wall conditions are $\Mach = 4.98, \Reyn=1477,  \Knud=0.005, T_{\infty}=184.6K, T_w = 1100K$. The hypersonic vehicle is placed in the centre of a spherical fluid domain with a diameter of $30 \, L$, where $L$ is the maximum length of the vehicle. In the computations, five meshes of decreasing element size are used, namely with 2105672, 5201671, 8519752 and 11007986 elements. In all computations the quadratic (P2) Dubiner basis functions are used. The finest mesh results in about 550 million degrees of freedom. For this example the boundaries have only been resolved linearly because our current high-order  boundary approximation does not allow for sharp features, such as corners  and edges. The computed pressure distribution on the vehicle surface and the cross-section of the flow are shown in Figure~\ref{fig:pressHypersonic} for the coarsest mesh with 2105672 elements. The pressure coefficient $C_p$ distribution along the symmetry axis on the bottom surface of the vehicle is plotted in Figure~\ref{fig:pressCoefHypersonic}. In the same figure, the experimentally obtained pressure distribution is also shown. The experiment and the computations exhibit an excellent  agreement.

\section{Conclusions }
%
We introduced a robust immersed discontinuous Galerkin method for solving the  compressible Navier-Stokes equations on unstructured non-boundary fitted meshes. Despite their many advantages, discontinuous Galerkin methods  suffer from their reliance on boundary-fitted meshes. This requirement has become a major bottleneck towards the realisation of automated design and simulation workflows in industrial applications with complex geometries and moving boundaries. As discussed throughout the paper, in discontinuous Galerkin methods there is no intrinsic need for boundary-fitted meshes. It is straightforward to apply the discontinuous Galerkin method on non-boundary fitted meshes by suitably adapting the evaluation of the element integrals and boundary fluxes in cut-elements. The stability or conditioning issues arising from cut-elements partially covered by the fluid can be effectively sidestepped with the introduced two cut-element stabilisation techniques. As our numerical studies confirm the proposed cut-element stabilisation techniques also work well in combination with the employed limiter for the fluxes. There are no spurious oscillations present close to the boundaries and there is no need to reduce the time step size in cut-elements.  In addition, the comparison of the computed flow parameters with experiments and previous computations shows good agreement. Amongst the two investigated cut-element stabilisation techniques, especially the scaling of the basis functions is straightforward to implement.

Our present high-order boundary reconstruction approach relies on implicitly defined smooth geometries. Most industrial CAD geometries consist however of trimmed NURBS (Non-Uniform Rational B-Spline) patches or other emerging geometry representations, like subdivision surfaces \cite{patrikalakis2009shape,Zorin:2000aa}. In addition, industrial geometries usually have many sharp features in form of corners and edges and contain many small geometry features, including chamfers, fillets and holes,  that cannot be realistically resolved with a discretisation mesh.  There is already extensive amount of work in computer graphics on the conversion between different geometry representations, the treatment of sharp features and multiresolution geometry representations, see e.g.~\cite{ Kobbelt:1998aa, kobbelt2001feature}. The application of these techniques in immersed discretisation methods aiming an achievable trade-off between accuracy and robustness is presently an active area of research~\cite{ruberg2016unstructured, bandara2016shape}.

\section*{Acknowledgement}

This study is partly supported by the Ministry of Education of the People's Republic of China as part of the China 111 project under grant number of B17037.

\bibliographystyle{wileyj}
\bibliography{immersedDG.bib}

\end{document}